\documentclass[aop,MSNbibl,seceqn,citesort,dvips]{arximspdf}
\usepackage{graphicx}

\doi{10.1214/11-AOP720} 
\volume{41}
\issue{3B}
\pubyear{2013}
\firstpage{2013}
\lastpage{2046}

\makeatletter

\newcommand{\calB}{\mathcal{B}}
\newcommand{\calC}{\mathcal{C}}
\newcommand{\calE}{\mathcal{E}}
\newcommand{\calF}{\mathcal{F}}
\newcommand{\calH}{\mathcal{H}}
\newcommand{\calL}{\mathcal{L}}
\newcommand{\calM}{\mathcal{M}}
\newcommand{\calN}{\mathcal{N}}
\newcommand{\calP}{\mathcal{P}}
\newcommand{\calR}{\mathcal{R}}
\newcommand{\calS}{\mathcal{S}}
\newcommand{\calT}{\mathcal{T}}
\newcommand{\calX}{\mathcal{X}}
\newcommand{\calY}{\mathcal{Y}}

\newcommand{\bbA}{\mathbb{A}}
\newcommand{\bbB}{\mathbb{B}}
\newcommand{\bbC}{\mathbb{C}}
\newcommand{\bbE}{\mathbb{E}}
\newcommand{\bbN}{\mathbb{N}}
\newcommand{\bbP}{\mathbb{P}}
\newcommand{\bbR}{\mathbb{R}}
\newcommand{\bbS}{\mathbb{S}}
\newcommand{\bbZ}{\mathbb{Z}}

\newcommand{\conn}{\leftrightarrow}

\newcommand{\dd}{\mathrm{d}}
\newcommand{\Piv}{\mathrm{Piv}}
\newcommand{\rwP}{\mathrm{P}}
\newcommand{\trwP}{{Q}}
\newcommand{\rwE}{\mathrm{E}}
\newcommand{\trwE}{\mathrm{E}_Q}
\newcommand{\rwZ}{\mathrm{Z}}
\newcommand{\Cconn}{\leftrightsquigarrow}

\newcommand{\unitvect}{\mathbf{e}_1}

\newtheorem{theorem}{Theorem}[section]
\newtheorem{lem}[theorem]{Lemma}
\newtheorem{prop}[theorem]{Proposition}
\newtheorem{cor}[theorem]{Corollary}

\newproclaim{rem}[theorem]{Remark}

\newcommand{\edges}{\mathsf{E}}

\makeatother

\begin{document}
\begin{frontmatter}

\title{Subcritical percolation with a line of defects}
\runtitle{Subcritical percolation with a line of defects}

\begin{aug}
\author[A]{\fnms{S.} \snm{Friedli}\thanksref{t1}\ead[label=e1]{sacha@mat.ufmg.br}},
\author[B]{\fnms{D.} \snm{Ioffe}\thanksref{t2}\ead[label=e2]{ieioffe@ie.technion.ac.il}}
\and
\author[C]{\fnms{Y.} \snm{Velenik}\corref{}\thanksref{t1}\ead[label=e3]{Yvan.Velenik@unige.ch}}
\runauthor{S. Friedli, D. Ioffe and Y. Velenik}
\affiliation{UFMG-ICEx, Technion and Universit\'e de Gen\`eve}
\address[A]{S. Friedli\\
Departamento de Matem\'atica\\
Avenida Ant\^{o}nio Carlos 6627\\
UFMG-ICEx C.P. 702\\
Belo Horizonte\\
30123-970 MG\\
Brasil\\
\printead{e1}}
\address[B]{D. Ioffe\\
Faculty of Industrial Engineering\\
\quad and Management\\
Technion, Haifa 32000\\
Israel\\
\printead{e2}}
\address[C]{Y. Velenik\\
Section de Math\'ematiques\\
Universit\'e de Gen\`eve\\
2-4 Rue du Li\`evre\\
1211 Gen\`eve 4\\
Switzerland\\
\printead{e3}} 
\end{aug}

\thankstext{t1}{Supported in part by the Swiss
National Science Foundation.}

\thankstext{t2}{Supported by the Israeli Science Foundation Grant
817/09.}

\received{\smonth{3} \syear{2011}}
\revised{\smonth{8} \syear{2011}}

\begin{abstract}
We consider the Bernoulli bond percolation process $\bbP_{p,p'}$ on the
nearest-neighbor edges of $\bbZ^d$, which are open independently with
probability $p<p_{c}$, except for those lying on the first
coordinate axis, for which this probability is $p'$. Define
\[
\xi_{p,p'}:=-\lim_{n\to\infty} n^{-1}\log\bbP_{p,p'}(0\conn n\unitvect)
\]
and $\xi_p:=\xi_{p,p}$. We show that there exists $p_c'=p_c'(p,d)$ such
that $\xi_{p,p'}=\xi_p$ if $p'<p_c'$ and $\xi_{p,p'}<\xi_p$ if
$p'>p_c'$. Moreover, $p_c'(p,2)=p_c'(p,3)=p$, and $p_c'(p,d)>p$ for
$d\geq 4$. We also analyze the behavior of $\xi_p-\xi_{p,p'}$ as
$p'\downarrow p_c'$ in dimensions $d=2,3$. Finally, we prove that when
$p'>p_c'$, the following purely exponential asymptotics holds:
\[
\bbP_{p,p'}(0\conn n\unitvect) = \psi_d  e^{-\xi_{p,p'} n}\bigl(1+o(1)\bigr)
\]
for some constant $\psi_d=\psi_d(p,p')$, uniformly for large values of
$n$. This work gives the first results on the rigorous analysis of
pinning-type problems, that go beyond the effective models and don't
rely on exact computations.
\end{abstract}

%

\end{frontmatter}

\section{Introduction and results}

We consider bond percolation on $\edges^d$, the set of nearest-neighbor
edges of $\bbZ^d$, $d\geq2$.
Let $\calL\subset\edges^ d$ be the set of all edges that lie on the
first coordinate axis $\{s\unitvect,s\in\bbR\}$, where $\unitvect$
denotes the unit vector $(1,0,\ldots,0)\in\bbR^d$.
Let $\bbP_{p,p'}$ be the probability measure on sets of configurations
of edges $\omega\in\{0,1\}^{\edges^d}$, under which each edge $e\in
\edges^d$ is open independently with probability
%
\begin{equation}\label{eqdefPppprime}
\bbP_{p,p'}\bigl(\omega(e)=1\bigr) =
\cases{
p, &\quad if $e\in\edges^ d\setminus\calL\equiv\calL^c$,\cr
p', &\quad if $e\in\calL$.}
\end{equation}

When $p'=p$, we write $\bbP_p$ instead of $\bbP_{p,p}$, and the model
coincides with ordinary homogeneous
Bernoulli edge percolation, whose critical threshold will be denoted
$p_c=p_c(d)$.

As far as we know, the properties of the connectivities under $\bbP
_{p,p'}$ were first studied by Zhang~\cite{Zhang}, who showed that in
$d=2$, there is no percolation under
$\bbP_{p_c(2),p'}$, for all $p'<1$.
Newman and Wu~\cite{NewmanWu} studied the same problem in large
dimensions as well as related properties, where the line $\calL$ is
replaced by higher-dimensional subspaces of $\bbZ^d$.

Let $\bbS^{d-1}$ be the unit sphere in $\bbR^d$.
It is well known~\cite{AB} that in the homogeneous case, for $p<p_c$,
\[
\xi_p(\mathbf{n}) := -\lim_{k\to\infty} \frac1k \log\bbP_p
(0\conn
[k\mathbf{n}])
\]
defines a function $\xi_p\dvtx\bbS^{d-1}\to(0,\infty)$ which can be
extended by positive homogeneity to a norm on $\bbR^d$. Let $\langle
\cdot,\cdot\rangle$ denote the inner product, and $|\cdot|$ the
Euclidean norm on $\bbR^d$.
There exists a convex, compact set $W_p\subset\bbR^d$ containing the
origin, such that for all $x\in\bbR^d$,
%
\begin{equation}\label{eqWulff}
\xi_p(x)=\sup_{t\in\partial W_p}\langle t,x\rangle.
\end{equation}
The sharp triangle inequality is also satisfied~\cite{CIV-Potts}: there
exists a constant
$c_1=c_1(p,d)>0$ such that
for all $x,y\in\bbR^d$,
%
\begin{equation}\label{eqSharpTriangle}
\xi_p(x)+\xi_p(y)-\xi_p(x+y)\geq c_1(|x|+|y|-|x+y|) .
\end{equation}

We also have, for any $x\in\bbZ^d$,
%
\begin{equation}\label{equpperboundxi}
\bbP_p(0\conn x)\leq e^{-\xi_p (x)} .
\end{equation}
It is also known~\cite{CampaninoIoffe} that the following
Ornstein--Zernike asymptotics holds, uniformly as $|x|\to\infty$:
%
\begin{equation}\label{eqOZsimple}
\bbP_p(0\conn x)=\frac{\Psi_d(x/|x|)}{|x|^{(d-1)/2}}e^{-\xi_p (x)}\bigl(1+o(1)\bigr),
\end{equation}
where $\Psi_d$ is a positive, real analytic function on $\bbS^{d-1}$.

Let ${\mathbf e}_j$, $j=1,\ldots,d$, denote the canonical basis of
$\bbR^d$.
By the symmetries of the lattice, $\xi_p(\unitvect)=\cdots=\xi
_p(\mathbf
{e}_d)$, and we define
%
\begin{equation}
\xi_p:= \xi_p(\unitvect) .
\end{equation}

In the inhomogeneous case, $p'\neq p$, the central quantity in
our analysis will be the modified inverse correlation length
%
\begin{equation}\label{eqxippp}
\xi_{p,p'} := -\lim_{n\to\infty} \frac1n\log\bbP_{p,p'}(0 \conn
n\unitvect) .
\end{equation}
Our goal is to study, for fixed $p<p_c$, the effect of the line $\calL$
on the rate of exponential decay $\xi_{p,p'}$. In particular, for which
%
\begin{figure}

\includegraphics{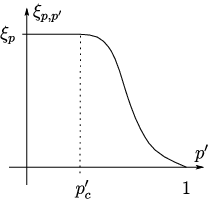}

\caption{A qualitative plot of $p'\mapsto\xi_{p,p'}$, for $d=2,3$.}
\label{figthm}
\end{figure}
values of $p'$ does $\xi_{p,p'}\neq\xi_p$? Our first main result is
the following; see also Figure~\ref{figthm}.
%
\begin{theorem}\label{thmmain}
Assume that $d\geq2$, $p<p_{c}$.
{\renewcommand\thelonglist{(\arabic{longlist})}
\renewcommand\labellonglist{\thelonglist}
\begin{longlist}
\item\label{mainitem1} The limit in (\ref{eqxippp}) exists for all
$0\leq p'\leq1$. Moreover, $p'\mapsto\xi_{p,p'}$ is Lipschitz
continuous and nonincreasing on $[0,1]$,
and
\[
\xi_{p,p'}>0\qquad \forall p'\in[0,1).
\]
\item\label{mainitem2} There exists $p'_c=p'_c(p,d)\in[p,1)$ such that
$\xi_{p,p'}=\xi_{p}$ for $p'\leq p'_c$ and $\xi_{p,p'}<\xi_{p}$ for
$p'>p'_c$.
On $(p_c',1)$, $p'\mapsto\xi_{p,p'}$ is real analytic and strictly decreasing.
\item\label{thmmainpoint1} When $d=2,3$, $p'_c=p$. Moreover, there
exist constants $c_2^\pm,c_3^\pm>0$ such that, as $p'\downarrow p_c'
= p$,
%
\begin{eqnarray}
\label{eqcriticald2}
c_2^-(p'-p)^2&\leq&\xi_p-\xi_{p,p'}\leq c_2^+(p'-
p
)^2 \qquad(d=2),\\
\label{eqcriticald3}
e^{-c_3^-/(p'-p
)}&\leq&\xi_p-\xi_{p,p'}\leq e^{-c_3^+/(p'- p)} \qquad(d=3).
\end{eqnarray}
\item\label{thmmainpoint2} When $d\geq4$, $p< p_c'<1$.
\end{longlist}}
\end{theorem}
%
\begin{rem}
Note that for $d=3$, (\ref{eqcriticald3}) rules out the possibility of
continuing $p'\mapsto\xi_{p,p'}$
analytically across $p$, to the interval $(0,p)$. It is an open
question whether such analytic continuation is possible in two dimensions.
\end{rem}
%
\begin{rem}\label{remnon-concavity}
We make a comment regarding the convexity/concavity of $p'\mapsto\xi
_{p,p'}$ for dimensions $2$ and $3$.
First, observe that $\xi_p$ diverges logarithmically as $p\downarrow
0$, and $\xi_{p,p'}\leq\xi_{0,p'}= |{\log p'}|$. Therefore, since in
dimensions $2$ and $3$ the slope of $\xi_{p,p'}$ (as a function of
$p'$) at $p_c'$ is equal to zero, there must be an inflection point
somewhere on the interval $(p_c',1)$, at least when $p$ is so small
that $\xi_p>1$. Note also that the above
implies that the Lipschitz constant must diverge at least as fast as
$|{\log p}|$, as $p\downarrow0$ (and at most as fast as $1/p$, as the
proof shows).
\end{rem}

In contrast to the polynomial correction in (\ref{eqOZsimple}) for the
homogeneous case, the presence of defects on the line $\calL$
leads to a purely exponential decay of the connectivities, which is the
content of our second result:
%
\begin{theorem}\label{thmmain2}
For all $d\geq2$ and for all $p'>p_c'$, there exists $\psi_d=\psi
_d(p,p')>0$ such that
%
\begin{equation}\label{eqpureexpon}
\bbP_{p,p'}(0\conn n\unitvect) = \psi_d e^{-\xi_{p,p'} n}
\bigl(1+o(1)\bigr) .
\end{equation}
\end{theorem}

As will be seen in Section~\ref{SecPureExpon}, the absence of a
polynomial correction in (\ref{eqpureexpon}) is due to the fact that
when $p'>p_c'$,
conditionally on $\{0\conn n\unitvect\}$, the cluster containing $0$
and $n\unitvect$, $C_{0,n\unitvect}$, is pinned
on the line $\calL$. Namely, as will be seen in Theorem~\ref
{propconeRenewal}, $C_{0,n\unitvect}$ splits into a string of
irreducible components centered on $\calL$ and whose sizes have
exponential tails.

The analysis of the effects of a line or a (hyper)plane of defects on
the qualitative statistical properties of polymers or interfaces has
been the subject of a large number of works dating back, at least, to
the late 1970s.
However, almost all rigorous studies to date have treated the framework
of effective models, in which the polymer/interface is modeled by the
trajectory of a random walk (or as a random function from $\bbZ^d\to
\bbR
$ in the case of higher-dimensional interfaces), and the understanding
of such models is by now very detailed~\cite{Giacomin,Velenik-PS}.
For example, in the case of a random walk pinned at the origin, one
studies the exponential divergence of the partition function
%
\begin{equation}\label{eqPartfuncRW}
Z_N^\varepsilon=E_{\mathit{RW}}[e^{\varepsilon L_N}|X_N=0] ,
\end{equation}
where $L_N$ is the local time of the random walk $X_k$ at the origin up
to time~$N$, and $\varepsilon>0$ is the pinning parameter (see Appendix
\ref{SecPinningRw}).

There is actually one very particular instance in which it has been
possible to investigate these phenomena in a noneffective setting: the
2d Ising model. Indeed, in this case it is sometimes possible to
compute explicitly the relevant quantities; see~\cite{Abraham} and
references therein. Needless to say, such computations do not convey
much understanding of the underlying physics (the desire to get a
better understanding of these exact results actually triggered the
analysis of effective models!).

On the other hand, new techniques developed during the last decade have
lead to a detailed description of structurally one-dimensional objects
in various lattice random fields, such as interfaces in 2d Ising and
Potts models~\cite{CIV-Ising,CIV-Potts,GI}, large subcritical clusters
in (FK-)percolation~\cite{CIV-Potts}, stretched self-interacting
polymers~\cite{IV-annealed}, etc.

The effect of a defect line in various systems has recently been the
focus of interest in different areas. In particular,
Beffara et al.~\cite{BSSS} have started to investigate the influence of
defects in the framework of last passage percolation.

It is worthwhile to point out an issue that makes the problem studied
in the present paper substantially more subtle than its effective
counterpart (\ref{eqPartfuncRW}). Namely, a natural way to compare
$\xi
_{p,p'}$ with $\xi_p$\vadjust{\goodbreak}
is to extract an effective weight for the cluster $C_{0,n\unitvect}$
connecting $0$ and $n\unitvect$. That is,
%
\begin{eqnarray}\label{eqmeasurerepulsive}
\frac{ \bbP_{p,p'} (0\conn n\unitvect) }{ \bbP_{p} (0\conn
n\unitvect) }
&=&\sum_{C\ni\{0,n\unitvect\}}\frac{\bbP_{p,p'}(C_{0,n\unitvect
}=C)}{\bbP_{p} (0\conn n\unitvect) }\nonumber\\
&=&\sum_{C\ni\{0,n\unitvect\}}e^{I(C)}\frac{\bbP
_{p}(C_{0,n\unitvect
}=C)}{\bbP_{p} (0\conn n\unitvect) }\\
&=& \bbE_p \bigl[ e^{I(C_{0,n\unitvect})} | 0\conn n\unitvect
\bigr],\nonumber
\end{eqnarray}
where
\[
I(C):= |C\cap\calL|\log\frac{p'}{p} + |\partial C\cap\calL|\log
\frac{1-p'}{1-p} ,
\]
and $\partial C$ denotes the exterior boundary of the cluster $C$, that
is, the set of all edges of $\edges^d\setminus C$ sharing at least one
endpoint with some edge of $C$.
Now, observe that in spite of the close resemblance of (\ref
{eqmeasurerepulsive}) with (\ref{eqPartfuncRW}), there is one major
difference: since
$\log\frac{p'}{p}$ and $\log\frac{1-p'}{1-p}$ always have opposite
signs, the effective interaction between the cluster and the line
$\calL
$ \textit{has both attractive and repulsive components}. This is a
manifestation of the presence of the ``phases'' that are neglected in
effective models, in which only the polymer/interface is considered and
not its environment.

Our analysis of
$\bbP_{p,p'}(0\conn n\unitvect)$ is based on the use of a geometrical
representation of the cluster $C_{0,n\unitvect}$ as an effective
directed random walk.
To use this representation effectively
for the lower bounds of part~\ref{mainitem2} of Theorem~\ref{thmmain},
the repulsive interaction of the cluster with $\calL$
will be handled with a suitable use of the Russo formula.

Random walk representations of subcritical clusters have been used in
\cite{CampaninoChayesChayes,CampaninoIoffe}
and~\cite{CIV-Potts}. The one used here is taken from~\cite{CIV-Potts},
and will be described in Section~\ref{SecRWREP}.
Standard renewal arguments are also recurrent in the paper; a reminder
of the main ideas can be found in Appendix~\ref{Apprenewal}.

\subsection{Open problems}
Although the picture provided by the present work is quite extensive,
we list here some open problems that we think would be particularly
interesting to investigate.

\begin{enumerate}[(P3)]
\item[(P1)] Properties of $\xi_{p,p'}$:
\begin{enumerate}[(a)]
\item[(a)] Analyze the behavior of $\xi_{p,p'}$ as $p'\downarrow
p'_{c}$, in dimensions $d\geq4$. In particular, determine whether
$\liminf_{p'\downarrow p_c'}\frac {\dd \xi_{p,p'}}{\dd p'}<0$ (which we
expect to be true in $d\geq6$, in analogy with the effective case
\cite{Giacomin}).
\item[(b)] Analyze the behavior of $\xi_{p,p'}$ as a function of both $p$
and $p'$. In particular, for $(p,p')$ close to the critical line
$p\mapsto p_c'(p)$.
\item[(c)] Determine, for all $p'\leq p'_{c}$, the sharp asymptotics of
the connectivity function $\bbP_{p,p'}(0\conn n\unitvect)$, and the
corresponding scaling limit of the cluster $C_{0,n\unitvect}$.\vadjust{\goodbreak}
\end{enumerate}

\item[(P2)] Introduce disorder, in which the occupation probabilities
of the edges $e\in\calL$ are i.i.d. random variables.
Study the relevance of disorder.
{For analogous considerations in the effective/directed case, we refer
to~\cite{AlexZyg,Giacomin,GLT} and references therein.}

\item[(P3)] More general defects:

\begin{enumerate}[(a)]
\item[(a)] Allow a defect line not coinciding with a coordinate axis, which
should be amenable to a rather straightforward adaptation of our techniques.
Or, as in~\cite{NewmanWu}, consider higher-dimensional defects like
hyperplanes of given codimension.
\item[(b)] Consider half-space percolation, with the defect line (or
hyperplane) at the boundary of the system. Although less natural from
the percolation point of view, such a setting is
relevant for the analysis of wetting phenomena.
\end{enumerate}

\item[(P4)] In each of the cases mentioned above, study the
connectivity $\bbP_{p,p'}(x\conn y)$ for generic points $x,y\in\bbZ^d$.

\item[(P5)] Extension to other models. In particular, a version for
FK-percolation seems feasible and would provide an extension of our
results to Ising/\break Potts models, which would be very interesting.
\end{enumerate}

We assume throughout the paper that edges outside $\calL$ are open with
probability~$p$, where $p<p_c$ is fixed.
Furthermore, $c_i$, $i=2,3,\ldots,$
will denote constants that can depend on the dimension $d$, on $p$ or
$p'$, but which remains uniformly bounded away from $0$ and $\infty$
for $(p,p')$ belonging to
compact subsets of $(0,p_c)\times(0,1)$.

{The line $\calL$ will often be identified with $\bbZ$. We will
therefore use the usual terminology related to the total order on $\bbZ
$ (such as ``being to the left of'' or ``being the largest among a set
of points'').
We will also consider $\calL$, without mention, sometimes as a set of
edges, and sometimes as a set of sites.}

\section{\texorpdfstring{Basic properties of $\xi_{p,p'}$}{Basic properties of xi p,p'}}

In this section, we prove items~\ref{mainitem1} and~\ref{mainitem2}
of Theorem~\ref{thmmain}, except for the strict monotonicity and
analyticity of $\xi_{p,p'}$, which will be proved, respectively, in
Sections~\ref{Ssecstrictmon} and~\ref{Ssecanalyticity}.

\textit{$\vartriangleright$ Existence of the limit.} The
existence of the limit in (\ref{eqxippp}) follows from the
subadditivity of the sequence
$n\mapsto-{\log\bbP_{p,p'}}(0\conn n\unitvect)$.

\textit{$\vartriangleright$ Monotonicity in $p'$ of $\xi
_{p,p'}$.} This follows from a standard coupling argument: if $p_1'\leq
p_2'$, then $\bbP_{p,p_1'}\preccurlyeq\bbP_{p,p_2'}$.

\textit{$\vartriangleright$ $\xi_{p,p'}=\xi_p$ for all
$p'\leq p$.}
Since $\xi_{p,p'}\geq\xi_p$ when $p'\leq p$, we only need to verify
that the reverse inequality also holds.
Let $0':=[n^\alpha] {\mathbf{e}_2}$ and $x':=n\unitvect+[n^\alpha]
{\mathbf{e}_2}$, where $1/2<\alpha<1$. We can realize $\{0\conn
n\unitvect\}$ by connecting $0$ to $0'$ and $n\unitvect$ to $x'$ by
straight segments of open edges, and by then
connecting $0'$ to $x'$ by an open path:
$\bbP_{p,p'}(0\conn n\unitvect)\geq(p^{n^\alpha})^2\bbP
_{p,p'}(0'\conn x')$.
If we characterize the event $\{0'\leftrightarrow x'\}$ by the
existence of a self-avoiding path $\pi\dvtx0'\to x'$,
\begin{eqnarray*}
\bbP_{p,p'}(0'\conn x')&\geq&\bbP_{p,p'}(\exists\pi\dvtx0'\to x', \pi
\cap
\calL=\varnothing)\\
&=&\bbP_{p}(\exists\pi\dvtx0'\to x', \pi\cap\calL=\varnothing).
\end{eqnarray*}
But by the van den Berg--Kesten (BK) inequality, (\ref{eqOZsimple})
and the {sharp triangle inequality (\ref{eqSharpTriangle})},
\begin{eqnarray*}
\bbP_{p}(\exists\pi\dvtx0'\leftrightarrow x', \pi\cap\calL\neq
\varnothing
)&\leq&\sum_{u\in\calL}\bbP_p(0\conn u)\bbP_p(u\conn x')\\
&\leq&\sum_{u\in\calL}e^{-\xi_p(u-0')-\xi_p(x'-u)}\\
&\leq& e^{-\xi_p(x'-0')}\sum_{u\in\calL}e^{-c_1(\|u-0'\|+\|x'-u\|-\|
x'-0'\|)}\\
&=&e^{-O(n^{2\alpha-1})}\bbP_p(0\conn n\unitvect) .
\end{eqnarray*}
Therefore, $\bbP_{p,p'}(0\conn n\unitvect)\geq p^{2n^\alpha
}(1-e^{-O(n^{2\alpha-1})})\bbP_p(0\conn n\unitvect)$,
which implies $\xi_{p,p'}\leq\xi_p$.

\textit{$\vartriangleright$ $\xi_{p,p'}<\xi_p$ for all $p'$
close enough to $1$.}
Namely, if $p'>e^{-\xi_p}$, then
by opening all the edges of $\calL$ between $0$ and $n\unitvect$,
\[
\bbP_{p,p'} (0\conn n\unitvect)
\geq{p'}^{n}= e^{(\log p'+\xi_p)n}e^{-\xi_pn}\geq e^{(\log p'+\xi_p)n}
\bbP_{p} (0\conn n\unitvect).
\]
The critical value
\[
p_c'=p_c'(p,d):=\sup\{p'\in[0,1]\dvtx\xi_{p,p'}=\xi_p\}
\]
thus separates the regime $\xi_{p,p'}=\xi_p$
from the one in which $\xi_{p,p'}<\xi_p$.

\textit{$\vartriangleright$ $\xi_{p,p'}>0$ for all $0\leq
p'<1$.} Define the slab
\[
\calS_{u,v}:=\{z\in\bbR^d \dvtx \langle u,\unitvect\rangle\leq
\langle
z,\unitvect\rangle<
\langle v,\unitvect\rangle\}.
\]
We divide $\calL_n:= \calL\cap\calS_{0,n\unitvect}$
into blocks of equal lengths $R\in\bbN$:
$B_j:=\calL_n\cap\calS_{jR\unitvect,(j+1)R\unitvect}$, with
$j=0,\ldots,[n/R]$.
Let also $\calH^-_j=\{x\dvtx\langle x,\unitvect\rangle< jR\}$,
$\calH^+_j=\{x\dvtx\langle x,\unitvect\rangle\geq(j+1)R\}$.
We say that $B_j$ is clear if there exists no path of open edges in
$\calL^c$ connecting $\calL\cap\calH_j^-$ to $\calL\cap\calH
_j^+$. We have
%
\begin{equation}\label{eqboundclear}\qquad
\bbP_{p,p'}(0\conn n\unitvect)\leq\bbP_{p,p'}(\mbox{each clear block
has at least one open edge}) .
\end{equation}
We show that when $R$ is large, a positive fraction of blocks is clear
with high probability.
For a cluster $C$ contained in $\calL^c$, let us define $l (C)$ and $r
(C)$ as, respectively, the left-most and right-most points of
intersections of the vertex
set of $C$ with $\calL$. We say that such $C$ is an $(R, n)$-bridge if
$r-l\geq R$ and the intersection $[l, r]\cap\calL_n\neq\varnothing$.
Let $C_1,C_2,\ldots, C_M$ be an enumeration of the disjoint
$(R, n)$-bridges. We set $l_i = l (C_i )$ and $r_i = r (C_i )$.
By construction, there are disjoint connections from $l_i$ to $r_i$ in\vadjust{\goodbreak}
$\calL^c; i=1, \ldots,M$.
If $0<\rho<1$, then, using the BK inequality in the last step,
\begin{eqnarray*}
&&\sum_{m=1}^\infty
\bbP_{p,p'}\Biggl( M= m ; \sum_{i=1}^m \bigl(r (C_i ) - l (C_i )\bigr) \geq\rho
n\Biggr)
\\
&&\qquad \leq
\sum_{m=1}^\infty\mathop{\mathop{\sum_{l_1,r_1,\ldots,l_m,r_m}}_{
r_i -l_i > R,}}_{\sum_{i=1}^m (r_i -l_i) \geq\rho n}
\bbP_{p,p'}\bigl(\circ_1^m \bigl\{ l_i\stackrel{\calL^c}{\conn}
r_i\bigr\} \bigr)
\\
&&\qquad \leq\sum_{m=1}^\infty
\mathop{\mathop{\sum_{l_1,r_1,\ldots,l_m,r_m}}_{r_i -l_i >
R,}}_{\sum
_{i=1}^m (r_i-l_i ) \geq\rho n}
\prod_{i=1}^m\bbP_{p,p'}\bigl(l_i\stackrel{\calL^c}{\conn} r_i\bigr) ,
\end{eqnarray*}
where it is understood that the points $l_i$ (resp., $r_i$) contributing
to the sum
are distinct, should in addition satisfy $[l_i ,r_i ]\cap\calL_n\neq
\varnothing$, and
$x\stackrel{A}{\conn}y$ means that $x$ and $y$ are connected by an open
path contained in $A$.
Now, $\bbP_{p,p'}(l_i\stackrel{\calL^c}{\conn} r_i)= \bbP
_{p}(l_i\stackrel{\calL^c}{\conn} r_i)\leq e^{-\xi_p|l_i-r_i|}$.
The contribution coming from segments so large that $[l_i ,r_i ]\supset
\calL_n$ is
clearly negligible, and we can restrict our attention to the case when
at least one of the
endpoints belongs to $\calL_n$.
Since, for all $t>0$, $\mathbf{1}_{\{{X\geq a}\}}\leq e^{t(X-a)}$,
this last sum is
bounded by
\begin{eqnarray*}
e^{-t\rho n}\sum_{m=1}^n\mathop{\sum_{l_1,r_1,\ldots
,l_m,r_m}}_{|l_i-r_i|> R}
\prod_{i=1}^me^{-(\xi_p-t)|l_i-r_i|}
&\leq& c_2e^{-t\rho n}\sum_{m=1}^n\sum_{l_1,\ldots,l_m}
2^m e^{-(\xi_p-t)mR}\\
&\leq& c_2e^{-t\rho n}\bigl(1+
2 e^{-(\xi_p-t)R}\bigr)^n
\end{eqnarray*}
for all $t<\xi_p$.\label{pagenummethod}
By taking $t=\xi_p/2$ and $R=\alpha/\xi_p$ with $\alpha$ large enough,
we get
\[
\bbP_{p,p'}\Biggl(\sum_{i=1}^M|l_i-r_i|\geq\rho n\Biggr)\leq c_2
e^{-\xi
_p\rho n/4} .
\]
This implies that
\[
\bbP_{p,p'}\bigl(\mbox{at least $[(1-\rho) n/2R]$ blocks are clear}\bigr)\geq
1-c_2 e^{-\xi_p\rho n/4} .
\]
Then, conditioned on the event that at least $[(1-\rho) n/2R]$ blocks
are clear, the probability on the right-hand side of (\ref
{eqboundclear}) is bounded above by\break
$\sum_{k\geq[(1-\rho) n/2R]}(1-(1-p')^R)^k\leq e^{-c_3 n}$. Altogether,
this shows that $\xi_{p,p'}>0$.

\textit{$\vartriangleright$ Lipschitz continuity of $\xi_{p,p'}$.}
The proof will rely on the following identity, which follows by Russo's
formula, and which will be used also later in Section~\ref{secLowerB}
(see~\cite{Grimmett}, page 44, for the proof of a similar
claim):

\begin{lem}\label{lemRusso}
For any increasing event $A$ with support in a finite subset of $\edges
^d$, and all $p_1', p_2'>0$,
%
\begin{equation}
\frac{\bbP_{p,p_2'}(A)}{\bbP_{p,p_1'}(A)}=\exp\int
_{p_1'}^{p_2'}\frac{1}{s}
\bbE_{p,s}[\# \Piv_{\calL}(A)|A] \,\dd s ,
\end{equation}
where $\Piv_{\calL}(A)$ is the set of pivotal edges $e\in\calL$ for
the event $A$.
\end{lem}

Let $\bbP_{p,p'}^{(n)}$ denote the restriction of $\bbP_{p,p'}$ to
the edges $\edges_n^d$ which lie in the
box $\Lambda_n:=[-a_n,a_n]^d\cap\bbZ^d$. {Since $\xi_{p,p'}>0$ for all
$p'<1$, we can
assume that $a_n\gg n$ is chosen sufficiently large so that for all $n$,}
%
\begin{equation}\label{eqfinitevolume}
\frac{1}{2}
\frac{\bbP_{p,p_2'}^{(n)}(0\conn n\unitvect)}{\bbP
_{p,p_1'}^{(n)}(0\conn n\unitvect)}\leq
\frac{\bbP_{p,p_2'}(0\conn n\unitvect)}{\bbP_{p,p_1'}(0\conn
n\unitvect
)}\leq2
\frac{\bbP_{p,p_2'}^{(n)}(0\conn n\unitvect)}{\bbP
_{p,p_1'}^{(n)}(0\conn n\unitvect)},
\end{equation}
when $n$ is large enough.
By Lemma~\ref{lemRusso}, for any $p_2'\geq p_1'\geq p_c'/2$,
%
\begin{equation}
\frac{\bbP_{p,p_2'}^{(n)}(0\conn n\unitvect)}{\bbP
_{p,p_1'}^{(n)}(0\conn n\unitvect)}=\exp\int_{p_1'}^{p_2'}\frac{1}{s}
\bbE_{p,s}^{(n)}[\# \Piv_{\calL}(0\conn n\unitvect)|0\conn
n\unitvect]
\,\dd s .
\end{equation}
Given a cluster $C_{0,n\unitvect}$, let $x$ (resp., $y$) be the leftmost
(resp., rightmost) site of $\calL\cap C_{0,n\unitvect}$, and
$L:=|x|+|y-n\unitvect|$.
We have
\begin{eqnarray*}
&&\bbE_{p,s}^{(n)}[\# \Piv_{\calL}(0\conn n\unitvect)|0\conn
n\unitvect
]\\
&&\qquad\leq
2n+e^{\xi_{p,s}(1+o(1))n}\sum_{\ell\geq n}(n+\ell)\bbP
_{p,s}^{(n)}(0\conn n\unitvect,L=\ell) .
\end{eqnarray*}
Since $\bbP_{p,s}^{(n)}(0\conn n\unitvect,L=\ell)\leq\ell\bbP
_{p,s}^{(n)}(0\conn(n+\ell)\unitvect)\leq\ell e^{-\xi_{p,s}(n+\ell
)}$, we get, using $p_c'\geq p$,
%
\begin{equation}
\frac{\bbP_{p,p_2'}^{(n)}(0\conn n\unitvect)}{\bbP
_{p,p_1'}^{(n)}(0\conn n\unitvect)}\leq\exp
\bigl(6(p_2'-p_1')n/p
\bigr) ,
\end{equation}
and thus $0\leq\xi_{p,p_1'}-\xi_{p,p_2'}\leq6(p_2'-p_1')/p$.

\section{Random walk representation of $C_{0,n\unitvect}$}\label{SecRWREP}

In this section we recall the description of $C_{0,n\unitvect}$ in
terms of a directed random walk, following~\cite{CIV-Potts}.
Since we only consider the direction $\unitvect$, the representation
simplifies in some respects. For instance, the inner products $\langle
y,t\rangle$ in~\cite{CIV-Potts} are replaced by $\langle y,\xi_p
\unitvect\rangle=\xi_p \langle y,\unitvect\rangle$.
The proofs of the main estimates {under $\bbP_{p}$} can be found in
\cite{CIV-Potts}.
The reader familiar with~\cite{CIV-Potts} can check the representation formulas
(\ref{eqrepresRW}), (\ref{eqrepresRWbis}) and (\ref{eqfundam}), and
proceed to Section~\ref{SecUpperbounds}.

{Observe that similar arguments for $\bbP_{p,p'}$ will be developed in
Section~\ref{SecPureExpon}.}\vadjust{\goodbreak}

\begin{figure}

\includegraphics{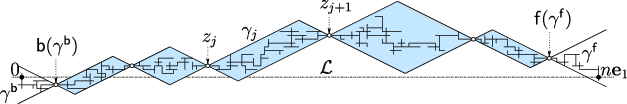}

\caption{The decomposition of $C_{0,n\unitvect}$ into irreducible components.}
\label{FigRWRepr}
\end{figure}

Let $0<\alpha<1$ be small enough so that the cone
%
\begin{equation}\label{eqdefCone}
\calY^>:=\{y\in\bbZ^d\dvtx \langle y,\xi_p\unitvect\rangle\geq
(1-\alpha
)\xi_p(y)\}
\end{equation}
has angular aperture at most $\pi/2$.
A point $z\in C_{0,n\unitvect}\neq\varnothing$ is called
cone-point if $0 <\langle z , \unitvect\rangle< n$ and
$C_{0,n\unitvect}\subseteq(z+ \calY^>)\cup(z- \calY^>)$.
We order the cone-points according to their first component:
$z_1,\ldots,z_{m+1}$.
By construction, $z_{i+1}\in z_{i}+\calY^>$.
The subgraphs
\[
\gamma_j:= C_{0,n\unitvect}\cap\calS_{z_{j}, z_{j+1}}
\]
are called cone-confined irreducible components of $C_{0,n\unitvect}$;
see Figure~\ref{FigRWRepr}. Note that $\gamma_j\subset
D(z_{j},z_{j+1})$, where
%
\begin{equation}
\label{eqdefDiammond}
D(z,z'):=(z+\calY^>)\cap(z'-\calY^>).
\end{equation}

The complement $C_{0,n\unitvect}\setminus(\gamma_1\cup\cdots\cup
\gamma
_m)$ can contain, at most, two connected components.
If it exists, the component containing $0$ (resp., $n\unitvect$) is
denoted $\gamma^{\mathsf b}$ (resp., $\gamma^{\mathsf f}$), and called backward
(resp., forward) irreducible.

Let ${\mathsf f}(\gamma_j):=z_j$ [resp., ${\mathsf
b}(\gamma_j):=z_{j+1}$]
denote the starting (resp., ending) point of $\gamma_j$, and ${\mathsf
f}(\gamma^{\mathsf f}):=z_m$, ${\mathsf b}(\gamma^{\mathsf b}):=z_1$.
Once a set of connected components $\gamma^{\mathrm b},\gamma_1,\ldots
,\gamma_m,\gamma^{\mathsf f}$ is given, compatible in the sense that
${\mathsf f}(\gamma_1)={\mathsf b}(\gamma^{\mathsf b})$, ${\mathsf b}(\gamma_m)={\mathsf
f}(\gamma^{\mathsf f})$, and
${\mathsf f}(\gamma_j)={\mathsf b}(\gamma_{j+1})$ if $j=1,\ldots,m-1$, then
these can be concatenated ($\sqcup$ denoting the corresponding
concatenation operation):
\[
\gamma^{\mathrm b}\sqcup\gamma_1\sqcup\cdots\sqcup\gamma_m\sqcup
\gamma^{\mathsf
f}\equiv C_{0,n\unitvect} .
\]
It can be shown that under $\bbP_p$, up to a term of order $e^{-\xi
_pn-\nu_1 n}$,
the number of cone-confined irreducible components grows linearly with $n$.

Therefore, the probability $\bbP_p(0\conn n\unitvect)$ can be
decomposed as
%
\begin{eqnarray}\label{eqrepresRW}
&&\bbP_p(0\conn n\unitvect) \nonumber\\[-8pt]\\[-8pt]
&&\qquad=\sum_{m\geq1}\mathop{\sum_{\gamma^{\mathrm b},\gamma_1,\ldots,\gamma
_m,\gamma^{\mathsf
f}}}_{\mathrm{compat.}}\bbP_p(C_{0,n\unitvect}=\gamma^{\mathrm b}\sqcup
\gamma_1\sqcup\cdots\sqcup\gamma_m\sqcup\gamma^{\mathsf f}) ,
\nonumber
\end{eqnarray}
where we neglected the configurations with less than two cone-points.
One can then define~\cite{CIV-Potts} independent events
$\Gamma^{\mathrm b}$, $\Gamma_1,\ldots,\Gamma_m$, $\Gamma^{\mathrm f}$ such that
%
\begin{equation}\label{eqrepresRWbis}\quad
\bbP_p(C_{0,n\unitvect}=\gamma^{\mathrm b}\sqcup\gamma_1\sqcup\cdots
\sqcup
\gamma_m\sqcup\gamma^{\mathsf f})
=\bbP_p(\Gamma^{\mathrm b})\Biggl(\prod_{j=1}^m\bbP_p(\Gamma_j)
\Biggr)\bbP
_p(\Gamma^{\mathrm f}) .
\end{equation}

The final step of the construction is to reformulate the rhs of (\ref
{eqrepresRW}) as the probability of an event involving a directed
random walk with independent increments. This follows a standard scheme
in renewal theory, sketched in Appendix~\ref{Apprenewal} in a simpler
situation, which starts by multiplying (\ref{eqrepresRW}) by $e^{\xi_pn}$.

First, we associate weights to the irreducible components $\gamma^{\mathsf
b}$ and $\gamma^{\mathsf f}$. By translation invariance, we can consider
$\gamma^{\mathsf f}$ as fixed at the origin,
and then translate it at~$n\unitvect$. If $u\in\calY^>$ and $v\in
-\calY^>$, define
%
\begin{equation}\label{eqweightsboundary}\qquad
\rho_{\mathsf b}(u):=e^{\langle u,\xi_p \unitvect\rangle}
\mathop{\sum_{\gamma^{\mathrm b}\ni0:}}_{{\mathsf b}(\gamma^{\mathsf
b})=u}\bbP
_p(\Gamma^{\mathsf b}) ,\qquad
\rho_{\mathsf f}(v):=e^{-\langle v,\xi_p \unitvect\rangle}
\mathop{\sum_{\gamma^{\mathrm f}\ni0:}}_{{\mathsf f}(\gamma^{\mathsf
f})=v}\bbP
_p(\Gamma^{\mathsf f}) .
\end{equation}
These weights satisfy
%
\begin{equation}\label{eqestimclusterbordv}
\rho_{\mathsf b}(u)\leq e^{-\nu_2|u|} ,\qquad
\rho_{\mathsf f}(v)\leq e^{-\nu_2|v|},
\end{equation}
where $\nu_2=\nu_2(p)>0$.
%
\begin{rem}\label{remboundaryterms}
Since the weights $\rho_{\mathsf b}$ and $\rho_{\mathsf f}$ have exponentially
decaying tails,
the sums over $u$ and $v$ [e.g., in the representation formulas
(\ref{eqfundam}) and (\ref{eqrepresRWbis}) below] can always
fix $\alpha>0$ small and
restrict attention to the terms for which $|u|,|v|\leq n^{1/2 - \alpha}$.
\end{rem}

Consider then the cone-confined components $\gamma_j$. Define the
displacement
\[
V(\gamma_j):={\mathsf b}(\gamma_j)-{\mathsf f}(\gamma_j) .
\]
By translation invariance, all components $\gamma_j$ with the same
displacement $V(\gamma_j)= y\in\calY^>$ have the same contribution to
the sum in (\ref{eqrepresRW}). We can thus consider only $\gamma_1$
and assume that its starting point is the origin: for all $y\in\calY^>$,
\[
{q}(y):=e^{\langle y,\xi_p \unitvect\rangle}
\mathop{\sum_{\gamma_1:}}_{{\mathsf f}(\gamma_1)=0, {\mathsf b}(\gamma
_1)=y}\bbP
_p(\Gamma_1)
.
\]
By a standard argument (a variant of Appendix~\ref{Apprenewal}), it
can be shown that $q$
defines a probability distribution on $\calY^>$.
Moreover, there exists $\nu_3=\nu_3(p)>0$ such that
%
\begin{equation}
\label{eqmassgap}
\sum_{y\dvtx|y|\geq t}q(y)\leq e^{-\nu_3 t} .
\end{equation}
Therefore, by summing over $u\in\calY^>$ and $v\in-\calY^>$, such
that $u_1<v_1$,
%
\begin{equation}\label{eqidentwithincrements}
e^{\xi_pn}\bbP_{p}(0\conn n\unitvect)=\\
\sum_{u,v}
\rho_{\mathsf b}(u)\rho_{\mathsf f}(v)
\sum_{m\geq1}
\mathop{\sum_{y_1,\ldots,y_m}}_{\sum_j
y_j=n\unitvect+v-u}
\prod_{j=1}^mq(y_j) .
\end{equation}
(As before, we neglected the term
with less than two cone-points.)\vadjust{\goodbreak}

Let us denote by $S=(S_k)_{k\geq0}$ the directed random walk on $\bbZ
^d$ whose
increments $Y_j=S_{j}-S_{j-1}\in\calY^>$ are i.i.d. and have
distribution $q$. When the walk is started at $u$, $S_0=u$, we denote
its distribution by $\rwP_u$. We can thus write (\ref
{eqidentwithincrements}) as
%
\begin{equation}\label{eqfundam}
e^{\xi_pn}\bbP_{p}(0\conn n\unitvect)=\sum_{u,v}\rho_{\mathsf
b}(u)\rho
_{\mathsf f}(u)\rwP_u\bigl(\calR(n\unitvect+v)\bigr),
\end{equation}
where
%
\begin{equation}
\calR(z):=\{\exists N\geq1 \mbox{ such that } S_N=z\} .
\end{equation}

More generally, if $A$ is an event measurable with respect to the
position of the endpoints of the irreducible components of
$C_{0,n\unitvect}$, that is, to the trajectory of the walk $S$, the
same procedure leads to
%
\begin{equation}\label{eqrepresRWbiss}
e^{\xi_pn}\bbP_{p}(A\cap\{0\conn n\unitvect\})
=\sum_{u,v}\rho_{\mathsf b}(u)\rho_{\mathsf f}(v)\rwP_u\bigl(A \cap\calR
(n\unitvect+v)\bigr).
\end{equation}

Let $Y_j=(Y_j^\parallel, Y_j^\perp)$ be the decomposition of $Y_j$ into
a longitudinal component $Y_j^\parallel:=\langle Y_j,\unitvect\rangle$
parallel to $\unitvect$, and a transverse component $Y_j^\perp\in
\bbZ
^{d-1}$ perpendicular to $\unitvect$. Then:
\begin{itemize}
\item$\rwP_u(Y_1^\parallel\geq1) = 1$;
\item$\rwP_u(|Y_1| > t) \leq e^{-\nu_3t}$ for large $t$;
\item for any $z^\bot\in\bbZ^{d-1}$, $\rwP_u(Y_1^\bot=z_\bot
) = \rwP_u(Y_1^\bot=-z_\bot)$.
\end{itemize}
Since the increments have exponential tails, the following local CLT
asymptotics along the direction $n\unitvect$ hold: fix $\alpha>0$.
Then, as $n\to\infty$,
%
\begin{equation}
\label{eqOZasympt}
\rwP_u\bigl(\calR(n\unitvect+v)\bigr)=
\frac{c_p}{n^{({d-1})/{2}}} \bigl(1+o(1)\bigr)
\end{equation}
for some constant $c_p>0$, uniformly in $|u|,|v|\leq n^{1/2 -
\alpha}$.
Together with (\ref{eqestimclusterbordv}) and (\ref{eqfundam}), this
in particular leads to the
Ornstein--Zernike asymptotics given in (\ref{eqOZsimple}) (for
$x=n\unitvect$).

\section{Upper bounds}\label{SecUpperbounds}
We now move on to the proof of the upper bounds of item \ref
{thmmainpoint1}, and of item~\ref{thmmainpoint2} of Theorem
\ref
{thmmain}.
We use (\ref{eqmeasurerepulsive}). Letting $\varepsilon:=\log(p'/p)>0$,
which is small if $p'-p$ is small, we get
%
\begin{equation}\label{equpperboundbasic}
\frac{\bbP_{p,p'}(0\conn n\mathbf{e}_1)}{\bbP_{p}(0\conn n\mathbf
{e}_1)} \leq\bbE_p\bigl[ e^{\varepsilon|C_{0,n\unitvect}\cap\calL|}
| 0\conn n\mathbf{e}_1 \bigr] .
\end{equation}
We use the random walk representation described in Section~\ref{SecRWREP}:
$C_{0,n\unitvect}=\gamma^{\mathsf b}\sqcup\gamma_1\sqcup\cdots\sqcup
\gamma
_m\sqcup\gamma^{\mathsf f}$.
If $S$ denotes the effective directed random walk associated to the
displacements of the components $\gamma_i$, we have
\begin{eqnarray*}
|C_{0,n\unitvect}\cap\calL|&=& |\gamma^{\mathsf b}\cap\calL|+|\gamma
^{\mathsf
f}\cap\calL|+\sum_{i=1}^m|\gamma_i\cap\calL|\\
&\leq&|\gamma^{\mathsf b}\cap\calL|+|\gamma^{\mathsf f}\cap\calL|+\sum
_{i=1}^m|D({S_{i-1},S_i})\cap\calL| ,
\end{eqnarray*}
{where the diamond $D(\cdot,\cdot)$ was defined in (\ref{eqdefDiammond}).}
If $\gamma^{\mathsf b}$ ends at $u$, define
$\rho_{\mathsf b}^\varepsilon(u)$
as in (\ref{eqweightsboundary}), with
$\bbP_p(\Gamma^{\mathsf b})e^{\varepsilon|\gamma^{\mathsf b}\cap\calL|}$
in place of $\bbP_p(\Gamma^{\mathsf b})$.
If $\gamma^{\mathsf f}$ starts at $v$,
$\rho_{\mathsf f}^\varepsilon(v)$ is defined in the same way. As can be verified,
exponential decay as in (\ref{eqestimclusterbordv}) holds for the weights
$\rho_{\mathsf f}^\varepsilon$ and $\rho_{\mathsf b}^\varepsilon$, when
$\varepsilon$ is sufficiently small. Still following Remark \ref
{remboundaryterms}, we will
only consider those $u,v$ with $|u|, |v|\leq n^{1/2 -\alpha}$
(for some $0<\alpha<1/2$).

Let $M:=\inf\{j\geq1\dvtx S_j=n\unitvect+v\}$.
Using (\ref{eqrepresRWbiss}), (\ref{eqOZasympt}) and (\ref{eqOZsimple}),
\[
\bbE_p\bigl[ e^{\varepsilon|C_{0,n\unitvect}\cap\calL|} |
0\conn
n\mathbf{e}_1 \bigr]
\leq c_4
\sum_{u,v}\rho_{\mathsf b}^\varepsilon(u)\rho_{\mathsf f}^\varepsilon(v)\rwE
_{u,v}\bigl[ e^{{\varepsilon\sum_{i=1}^M} |{D(S_{i-1},S_i) \cap\calL
}|}\bigr] ,
\]
where $\rwE_{u,v} [\cdot]=\rwE_u [\cdot|\calR(n\unitvect
+v)]$.
As we said,
%
\begin{equation}
\sum_{u,v}\rho_{\mathsf b}^\varepsilon(u)\rho_{\mathsf f}^\varepsilon(v)<\infty
.
\end{equation}
We further decompose
\[
\rwE_{u,v}\bigl[ e^{{\varepsilon\sum_{i=1}^M} |{D(S_{i-1},S_i) \cap
\calL
}|}\bigr]
= \sum_{m=1}^n \rwE_{u,v}\bigl[ e^{{\varepsilon\sum_{i=1}^m}
|{D(S_{i-1},S_i) \cap\calL}|},{M=m}\bigr] .
\]
Therefore, for all fixed $1\leq m_0\leq n$,
%
\begin{equation}\label{equpperboundEuv}\quad
\rwE_{u,v}\bigl[ e^{{\varepsilon\sum_{i=1}^M }|{D(S_{i-1},S_i) \cap
\calL
}|}\bigr]
\leq e^{\varepsilon n}\rwP_{u,v}(M< m_0)+n\sup_{m_0\leq m\leq
n}A_{u,v}(m) ,
\end{equation}
where
%
\begin{eqnarray}
A_{u,v}(m):\!&=&\rwE_{u,v}\bigl[ e^{{\varepsilon\sum_{i=1}^m} |{D(S_{i-1},S_i)
\cap\calL}|}\bigr]\nonumber\\
&=&\sum_{k=1}^m\mathop{\sum_{\ell_1,\ldots,\ell_k}}_{\sum_j\ell_j=m}
\rwE_{u,v}\Biggl[\prod_{j=1}^k\Psi_\calL(S_{a_{j}-1},S_{a_j})\Biggr]
\end{eqnarray}
with $\Psi_\calL(S_{i-1},S_{i}):=e^{\varepsilon|D(S_{i-1},S_i)\cap
\calL
|}-1$, and
where $\ell_j\geq1$, $a_j:=\ell_1+\cdots+\ell_j$, {$a_0:=0$}.
%
\begin{figure}

\includegraphics{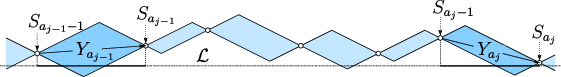}

\caption{The proof of the upper bound: the size of the intersection of
a diamond with
$\calL$ is bounded above by the size of its projection on $\calL$.}
\label{FigPinnedPoints}
\end{figure}
Remembering that the cone $\calY^>$ has an opening angle of at most
$\pi
/2$, we have (see Figure~\ref{FigPinnedPoints})
%
\begin{equation}\label{eqintersdiammond}
|{D(S_{i-1},S_i) \cap\calL}|\leq Y_i^\parallel\mathbf{1}_{\{
{|S_{i-1}^\bot| \leq{Y_i^\parallel}}\}} .
\end{equation}
Therefore,
\begin{eqnarray*}
\psi_\calL(S_{i-1},S_i)&\leq& e^{\varepsilon Y_i^\parallel\mathbf{1}_{\{
{|S_{i-1}^\bot| \leq{Y_i^\parallel}}\}}}-1\\
&=&
(e^{\varepsilon Y_i^\parallel}-1 ) \mathbf{1}_{\{{|{S_{i-1}^\bot}|
\leq Y_i^\parallel}\}}\\
&\leq& (e^\varepsilon-1)Y_i^\parallel e^{\varepsilon Y_i^\parallel}\mathbf
{1}_{\{{|{S_{i-1}^\bot}| \leq Y_i^\parallel}\}}\\
&\equiv& (e^\varepsilon-1)B(S_{i-1}^{\bot},Y_i) ,
\end{eqnarray*}
which yields
\begin{eqnarray*}
A_{u,v}(m)
&\leq& \sum_{k=1}^m(e^{\varepsilon}-1)^{k}\mathop{\sum_{\ell_1,\ldots
,\ell
_k}}_{\sum_j\ell_j=m}
\rwE_{u,v}\Biggl[\prod_{j=1}^kB(S_{a_{j}-1}^{\bot},Y_{a_j})\Biggr]\\
&\leq& O\bigl(n^{({d-1})/{2}}\bigr)
\sum_{k=1}^m(e^{\varepsilon}-1)^{k}\mathop{\sum_{\ell_1,\ldots,\ell
_k}}_{\sum_j\ell_j=m}
\rwE_{u}\Biggl[\prod_{j=1}^kB(S_{a_{j}-1}^{\bot},Y_{a_j})\Biggr] ,
\end{eqnarray*}
where (\ref{eqOZasympt}) was used again.
For all $j$, by {the Markov property} and the local limit theorem in
dimension $d-1$ (see Figure~\ref{FigPinnedPoints} and note that the
upper bound below is trivial
whenever $a_{j-1} = a_j -1$),
\begin{eqnarray*}
\rwE_{u}[B(S_{a_{j}-1}^{\bot},Y_{a_j})|S_{a_{j-1}},Y_{a_j}]&=&
Y_{a_j}^\parallel e^{\varepsilon Y_{a_j}^\parallel}
\rwP_{S_{a_{j-1}}}(|{S_{a_{j}-1}^\bot}| \leq Y_{a_j}^\parallel)\\
&\leq& Y_{a_j}^\parallel e^{\varepsilon Y_{a_j}^\parallel}
c_5(Y_{a_j}^\parallel)^{d-1}\ell_j^{-({d-1})/{2}} .
\end{eqnarray*}
Therefore, since $\rwP_u(Y_{a_j}^\parallel\geq t)\leq e^{-\nu_3 t}$,
\[
\rwE_{u}[B(S_{a_{j}-1}^{\bot},Y_{a_j})|S_{a_{j-1}}]\leq
c_5\sum_{t\geq1}t^{d} e^{\varepsilon t}e^{-\nu_3 t}
\ell_j^{-({d-1})/{2}}\equiv c_6\ell_j^{-({d-1})/{2}}
\]
with $c_6<\infty$ if $\varepsilon<\nu_3$.
This gives $A_{u,v}(m)\leq O(n^{({d-1})/{2}})A(m)$, where
%
\begin{equation}\label{eqdefAm}
A(m):=
\sum_{k=1}^m\bigl(c_6(e^{\varepsilon}-1)\bigr)^{k}
\mathop{\sum_{\ell_1,\ldots,\ell_k}}_{\sum_j\ell_j=m}\prod
_{j=1}^k\ell
_j^{-({d-1})/{2}} .
\end{equation}
In dimensions $d\geq4$, we ignore the constraint $\sum_j\ell_j=m$ and
bound $A(m)$
uniformly by
\[
A(m)\leq\sum_{k=1}^\infty\biggl\{c_6(e^{\varepsilon}-1)
\sum_{\ell\geq1}\ell^{-({d-1})/{2}}
\biggr\}^{k} ,
\]
which converges when $\varepsilon>0$ is small enough. Therefore, using
(\ref{equpperboundEuv}) with $m_0=1$, (\ref{equpperboundbasic}) is
\[
\bbE_p\bigl[ e^{\varepsilon|C_{0,n\unitvect}\cap\calL|} |
0\conn
n\mathbf{e}_1 \bigr]
= O\bigl(n^{({d+1})/{2}}\bigr).
\]
This shows
that $\xi_{p,p'}\geq\xi_p$ when $p'-p$ is small enough. Combined with
$\xi_{p,p'}\leq\xi_p$, this implies that $p_c'(p,d)>p$ in dimensions
$d\geq4$.

In dimensions $d=2$ and $3$, we obtain an upper bound on $A(m)$ which
diverges with $m$, in a standard way. As in Appendix~\ref{Apprenewal},
consider the generating function
\[
\bbA(s):=\sum_{m\geq1}A(m)s^m .
\]
Using (\ref{eqdefAm}),
$\bbA(s)=\sum_{k\geq1}\bbB(s)^k$
where $\bbB(s):=c_6(e^{\varepsilon}-1)\sum_{\ell\geq1}\ell^{-
({d-1})/{2}}s^\ell$.
Let $\phi(\varepsilon)>0$ be the unique number for which
%
\begin{equation}\label{eqdefrepsilon}
\bbB\bigl(e^{-\phi(\varepsilon)}\bigr)=1 .
\end{equation}
We have $\bbA(e^{-2\phi(\varepsilon)})<\infty$, and therefore
$A(m)\leq
e^{2\phi(\varepsilon)m}$ for all large enough $m$. Using
(\ref{equpperboundEuv}) with $m_0=c_7 n$ with
$c_7>0$ small enough, and taking $\varepsilon$ small enough,
(\ref{equpperboundbasic}) is bounded by
\[
\bbE_p\bigl[ e^{\varepsilon|C_{0,n\unitvect}\cap\calL|} |
0\conn
n\mathbf{e}_1 \bigr]
\leq c_8\bigl(1+O\bigl(n^{({d-1})/{2}}\bigr)e^{2\phi(\varepsilon)n}\bigr).
\]
This shows that $\xi_p-\xi_{p,p'}\leq2\phi(\varepsilon)$.
Using~\cite{Giacomin}, Theorem A.2, in (\ref{eqdefrepsilon}), the
asymptotics of $\phi(\varepsilon)$ when $\varepsilon\downarrow0$ is seen
to be
\[
\phi(\varepsilon) =
\cases{
c_9\varepsilon^2 \bigl(1+o(1)\bigr) &\quad$(d=2)$,\vspace*{2pt}\cr
{e^{-c_{10}/\varepsilon(1+o(1))}} &\quad$(d=3)$.}
\]

\section{Lower bounds}\label{secLowerB}

We prove the lower bounds of item~\ref{thmmainpoint1} of Theorem
\ref{thmmain},
in $d=2,3$, for $p'>p$, with $p'-p$ small enough. We will need the
following rough estimate on the connectivity under $\bbP_{p,p'}$:
%
\begin{lem}\label{lemexpdecay}
Set
%
\begin{equation}
\label{eqdefxietoile}
\xi^*_p:=\min_{\mathbf{n}\in\bbS^{d-1}}\xi_p(\mathbf{n}) > 0 .
\end{equation}
For all $p<p_c$, there exists {$\eta=\eta(p)>0$} such that, for all
$p'<p+\eta$,
\[
\bbP_{p,p'}(x\conn y) \leq e^{-c_{11}\xi_p^*|y-x|},
\]
uniformly in $x,y\in\bbZ^d$.
\end{lem}
\begin{pf}
Let $\ell(x ,y ):=|C_{x ,y }\cap\calL|$. Proceeding as in (\ref
{eqmeasurerepulsive}),
\[
\bbP_{p,p'}(x\conn y ) \leq\bbE_p \bigl[ e^{\ell(x ,y )\log(p'/p) } ;
x\conn y \bigr].
\]
Since
$\bbP_{p}(x\conn y ; \ell(x ,y ) = l ) \leq e^{-c_{12} \xi_p^* |x-y
|\wedge l }$,
the claim follows.
\end{pf}

Recall that $\bbP_\cdot^{(n)}$ denotes the restriction of $\bbP
_\cdot$ to
the edges $\edges_n^d$ which lie inside a large box
$\Lambda_n$, so that by (\ref{eqfinitevolume})
\[
\frac{ \bbP_{p,p'} (0\conn n\unitvect) }{ \bbP_{p} (0\conn
n\unitvect)
} \geq\frac{1}{2}
\frac{ \bbP_{p,p'}^{(n)} (0\conn n\unitvect) }{ \bbP_{p}^{(n)}
(0\conn
n\unitvect) } .
\]

Let $\calP_n$ denote the collection of self-avoiding nearest-neighbor
paths $\pi\dvtx0\to n\unitvect$ contained in $\Lambda_n$.
Let $\pi=(\pi_0,\pi_1,\ldots,\pi_m)\in\calP_n$, that is,
$\pi_0=0$ and $\pi_m=n\unitvect$.
We say that $\pi_i$ is a cone-point of $\pi$ if
$0< \langle\pi_i ,\unitvect\rangle<n $
and $\pi\subset(\pi_i - \calY^>)\cup( \pi_i + \calY^> )$.

Let $\delta>0$, and define
\[
\calM_\delta:= \{\mbox{$\exists$ an open path $\pi\in\calP
_n$ with
at least $\delta n$
cone-points on $\calL_n$}\}.
\]
We emphasize the crucial fact that we do not require that cone-points
of open paths are
cone-points of the whole cluster $C_{0 ,n\unitvect}$. This ensures that
$\calM_\delta$
is an increasing event: once a configuration contains a suitable open
path, opening additional edges will never remove this path (observe
also that suitability of an open path only depends on its geometry, not
on the state of other edges in the configuration).

Since $ \{0\conn n\unitvect\}\supset\calM_\delta$, we can write
\[
\frac{ \bbP_{p,p'}^{(n)} (0\conn n\unitvect) }{ \bbP_{p}^{(n)}
(0\conn
n\unitvect) } \geq
\frac{ \bbP_{p,p'}^{(n)} (\calM_\delta)}{ \bbP_{p}^{(n)} (0\conn
n\unitvect) } =\frac{\bbP_{p,p'}^{(n)}(\calM_\delta)}{\bbP
_p^{(n)}(\calM
_\delta)}
\bbP_p^{(n)}(\calM_\delta| 0\conn n\unitvect).
\]
The terms in the last display are, respectively, the energy gain and
the entropy cost for restricting to the event $\calM_\delta$. These
will be studied separately. First:
%
\begin{prop}\label{PropGainBound} Let $d\geq2$.
There exists $c_{13}=c_{13}(p,p')>0$ such that, for all $p'>p$, $p'-p$
small enough, and all $n\in\bbN$,
\[
\frac{\bbP_{p,p'}^{(n)}(\calM_\delta)}{\bbP_p^{(n)}(\calM_\delta
)} \geq
e^{c_{13}\delta(p'-p)n}.
\]
\end{prop}

Then, we check that $\calM_\delta$ is not too unlikely under $\bbP
_p^{(n)}(\cdot|0\conn n\unitvect)$:
%
\begin{prop}\label{PropEntropyCost}
There exist $c_{14}=c_{14}(p)>0$ and $c_{15}=c_{15}(p)>0$ such that for
small enough $\delta>0$,
\[
\bbP_p^{(n)}(\calM_\delta| 0\conn n\unitvect) \geq
\cases{
e^{-c_{14} \delta^2 n} &\quad$(d=2)$,\cr
e^{-c_{15} (\delta/|{\log\delta}|) n} &\quad$(d=3)$.}\vadjust{\goodbreak}
\]
\end{prop}

Putting these bounds together, an appropriate choice of $\delta$ as a
function of $p'-p$
leads to the lower bounds of item~\ref{thmmainpoint1} of
Theorem~\ref{thmmain}. Namely,
\[
\delta:=
\cases{
c_{13}(p'-p)/(2c_{14}),&\quad in $d=2$,\vspace*{2pt}\cr
e^{-2c_{15}/(c_{13}(p'-p))},&\quad in $d=3$.}
\]

\subsection{\texorpdfstring{Proof of Proposition \protect\ref{PropGainBound}}{Proof of Proposition 5.2}}
\label{Subsecpaths}
First, observe that
\[
\Piv_{\calL\cap\Lambda_n}(\calM_\delta)
\supset
\Piv_{\calL\cap\Lambda_n}(0\conn n\unitvect)\qquad\mbox{on the event }\calM
_\delta.
\]
Indeed, let $e\in\Piv_{\calL\cap\Lambda_n}(0\conn n\unitvect )$. Then
$e$ must belong to \textit{all} paths $\pi$ satisfying the conditions
prescribed in the event $\calM_\delta$ (since removing this edge
disconnects $0$ from $n\unitvect$). This shows that $e$ is pivotal for
$\calM_\delta$.

We start by using Lemma~\ref{lemRusso}: by the preceding observation
and the fact that $\calM_\delta$ is increasing, we obtain
%
\begin{eqnarray}\label{eqalmostthere}
\frac{\bbP^{(n)}_{p,p'}(\calM_\delta)}{\bbP^{(n)}_p(\calM_\delta)}
&=&\exp\int_p^{p'}\frac1s \bbE^{(n)}_{p,s} [\#\Piv_{\calL\cap
\Lambda_n}(\calM_\delta) | \calM_\delta] \,\dd s \nonumber\\
&\geq&\exp\int_p^{p'}\frac1s \bbE^{(n)}_{p,s} [\#\Piv_{\calL
\cap
\Lambda_n}(0\conn n\unitvect) | \calM_\delta] \,\dd s
\\
&\geq&\exp\int_p^{p'}\frac1s \bbE^{(n)}_{p,s} [\#\Piv_{\calL
_n}(0\conn n\unitvect) | \calM_\delta] \,\dd s .\nonumber
\end{eqnarray}
Our goal is thus to bound $\#\Piv_{\calL_n}(0\conn n\unitvect)$ from
below on $\calM_\delta$.

Let us fix an arbitrary total ordering on $\calP_n$. For
each $\pi\in\calP_n$, let
$\calE_\pi\subset\calM_\delta$ denote the event on which $\pi$ is the
smallest open
path having at least $\delta n$ cone-points on $\calL_n$.
Then
%
\begin{eqnarray}\label{eqbigexpression}
&&\bbE^{(n)}_{p,s} [\#\Piv
_{\calL_n}(0\conn n\unitvect) | \calM_\delta]\nonumber\\[-8pt]\\[-8pt]
&&\qquad=\sum_{\pi\in\calP_n}
\bbE^{(n)}_{p,s} [\#\Piv_{\calL_n}(0\conn n\unitvect) |
\calE
_\pi]
\bbP^{(n)}_{p,s} (\calE_\pi|\calM_\delta) .
\nonumber
\end{eqnarray}
Let $\edges^d_{n,\pi}:=\edges^d_n\setminus\pi$.
We say that a cone-point $\pi_s\in\calL_n$ is covered if
\[
(\pi_s-\calY^>) \stackrel{\edges^d_{n,\pi}}{\longleftrightarrow}
(\pi
_{s}+\calY^>) ,
\]
uncovered otherwise.
%
\begin{lem}
\label{lemcovering}
Given $\pi\in\calP_n$ and $\rho>0$ define the event
\[
{\mathrm A}_\pi(\rho,n )=
\{
\mbox{A fraction $\geq\rho$ of $\pi$'s cone-points on $\calL_n$
are uncovered}\}.
\]
Let $s-p>0$ be sufficiently small.
Then there exists $\rho=\rho(p)>0$ such that
for all $\pi\in\calP_n$ {compatible with $\calM_\delta$},
%
\begin{equation}\label{eqlowerprobstretches}
\bbP_{p,s}^{(n)}\bigl({\mathrm A}_\pi(\rho,n )
|\calE_\pi\bigr) \geq\tfrac12.
\end{equation}
\end{lem}

Observe that each uncovered cone-point of $\pi$ on $\calL_n$ has two
incident edges $e\in\calL_n$ which are pivotal for $\{0\conn
n\unitvect
\}$. Therefore, by (\ref{eqlowerprobstretches}),
\[
\bbE^{(n)}_{p,s} [\#\Piv_{\calL_n}(0\conn n\unitvect) |
\calE
_\pi]\geq\tfrac{1}{2} \times\rho{\delta n} ,
\]
which, together with (\ref{eqalmostthere}) and
(\ref{eqbigexpression}), completes the proof of Proposition
\ref{PropGainBound}.
\begin{pf*}{Proof of Lemma~\ref{lemcovering}}
Fix some path $\pi$ realizing $\calM_\delta$.
We claim first that, as probability measures on $\{0,1\}^{\edges
^d_{n,\pi}}$,
%
\begin{equation}\label{eqcoupling}
\bbP_{p,s}^{(n)}(\cdot| \calE_\pi) \preccurlyeq\bbP
_{p,s}^{(n)}(\cdot) .
\end{equation}
Indeed, note that if $\omega\in\calE_\pi\subseteq\{
0,1\}
^{\edges^d_{n ,\pi}}$, then for every edge $e\in\edges_{n,\pi}^d$,
the configuration $\omega_e^0$ defined by
\[
\omega_e^0 (b) =
\cases{
\omega(b) , &\quad if $b\neq e$,\cr
0, &\quad if $b=e$,}
\]
belongs to $\calE_\pi$ as well. In particular, any two configurations
$\omega,\omega'\in\calE_\pi$ are
connected via a sequence of bond flips within $\calE_\pi$. Furthermore,
for every $\eta\in\calE_\pi$ and
for any edge $e\in\edges_{n, \pi}^d$,
\begin{eqnarray*}
&&\bbP_{p,s}^{(n)}\bigl(\omega(e)=1 |
\omega|_{\edges_{n ,\pi}^d\setminus\{e\}} = \eta;
\calE_\pi\bigr) \\
&&\qquad=
\cases{
0, &\quad if $e$ is pivotal for $\calE_\pi$ in $\eta$,\cr
\bbP_{p,s}^{(n)} \bigl(\omega(e)=1\bigr), &\quad otherwise.}
\end{eqnarray*}
Thus, (\ref{eqcoupling}) follows from a a standard dynamic coupling
argument for two
Markov chains on $\{ 0, 1\}^{\edges_{n ,\pi}^d}$, which are
reversible with respect to $\bbP_{p,s}^{(n)}(\cdot)$ and
$\bbP_{p,s}^{(n)} (\cdot|\calE_\pi)$ accordingly.

The event ${\mathrm A}_\pi(\rho, n )$ is $\edges_{n, \pi}^d$-measurable
and decreasing.
Hence, in order to prove (\ref{eqlowerprobstretches}) it would be
enough to
show that $\bbP_{p,s}^{(n)} ({\mathrm A}_\pi(\rho, n ))\geq\frac12$
for all $\calM_\delta$-compatible paths $\pi\in\calP_n$.

%
\begin{figure}[b]

\includegraphics{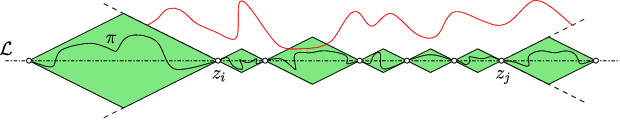}

\caption{The event $\{z_i\Cconn z_{j}\}$.}
\label{Figwiggle}
\end{figure}

Let us fix such a $\pi$, and denote the cone-points of $\pi$ on
$\calL
_n$, ordered from left to right, by
$z_1,\ldots,z_M$, $M\geq\delta n$.
We denote by
$z_i\Cconn z_j$ ($i<j$) the event (see Figure~\ref{Figwiggle})
\[
(z_i-\calY^>) \stackrel{\edges_{n,\pi}^d}{\longleftrightarrow}
(z_j+\calY^>) .
\]
By construction the events $z_i\Cconn z_j$ are $\edges_{n, \pi}^d$-measurable
and increasing.\vadjust{\goodbreak}

Observe that if $\pi$ has $m$ of its points $z_j$ covered,
then there must exist a set of distinct pairs $\{z_{k_j},z_{k_j'}\}
\subset\{z_1,\ldots,z_M\}$, $j=1,\ldots,q$, such that:
\begin{longlist}[(2)]
\item[(1)] $\sum_{j=1}^q|k_j'-k_j+1|=m$;
\item[(2)] $\{z_{k_1}\Cconn z_{k_1'}\}\circ\cdots\circ\{z_{k_q}\Cconn
z_{k_q'}\}$.
\end{longlist}
By the BK inequality,
\[
\bbP_{p,s}^{(n)} (\{z_{k_1}\Cconn z_{k_1'}\}\circ\cdots\circ\{
z_{k_q}\Cconn z_{k_q'}\})
\leq\prod_{j=1}^q \bbP_{p,s}^{(n)} (z_{k_j}\Cconn z_{k_j'}).
\]
Now, it follows from Lemma~\ref{lemexpdecay} that if $s-p$ is small enough,
and $|z_{k_j}- z_{k_j'}|\geq c_{16}/\xi_p^*$,
\begin{eqnarray*}
\bbP_{p,s}^{(n)} (z_{k_j}\Cconn z_{k_j'})
&\leq&\mathop{\sum_{x\in z_{k_j}-\calY^>}}_{y\in z_{k_j'}+\calY^>}
e^{-(c_{11}\xi_p^*/2)(|x-y|+1)}
\leq c_{17} e^{-(c_{11}\xi_p^*/2)|z_{k_j}- z_{k_j'}|}\\
&\leq& e^{-c_{18}(|{k_j'}-{k_j}|+1)} .
\end{eqnarray*}
On the other hand, if $|z_{k_j}- z_{k_j'}|\leq c_{16}/\xi_p^*$, then
\[
\bbP_{p,s}^{(n)} (z_{k_j}\Cconn z_{k_j'})\leq\bbP_{p,s}^{(n)}
(z_{k_j}\mbox{ is covered})\leq e^{-c_{19}}\leq
e^{-c_{20}(|{k_j'}-{k_j}|+1)} .
\]
Indeed, if
$B_R(z)$ is the Euclidean ball of radius $R$ centered at $z$, and
$\calB:=\{\mbox{all edges of $B_R(z_{k_j})$ are closed}\}$ with
$R=c_{21}/\xi_p^*$ with $c_{21}>0$ large enough, then
\[
\bbP_{p,s}^{(n)} (z_{k_j}\mbox{ is not covered})\geq
\bbP_{p,s}^{(n)} (z_{k_j}\mbox{ is covered}\, | \calB)\bbP
_{p,s}^{(n)}(\calB)\geq\tfrac{1}{2}\times e^{-c_{22}R^{d}} .
\]
Therefore, it follows from (\ref{eqcoupling}) and the above
discussion that
with $c_{23}:=c_{18}\wedge c_{20}$,
\begin{eqnarray*}
&&\bbP_{p,s}^{(n)} (\mbox{a fraction $\geq\alpha$ of $\pi$'s
cone-points on $\calL_n$ are covered}\,|\calE_\pi)\\
&&\qquad\leq\sum_{m=\alpha M}^M \pmatrix{M\cr m} e ^{-c_{23}m}\leq
e^{-c_{24} M}
\leq e^{-c_{24} \delta n}
\end{eqnarray*}
once $\alpha$ is close enough to $1$. This proves the lemma.
\end{pf*}

\subsection{\texorpdfstring{Proof of Proposition \protect\ref{PropEntropyCost}}{Proof of Proposition 5.3}}
\label{secRWrepresentation}
We use the representation of $C_{0,n\unitvect}$ in terms of
the directed random walk $S$. Observe that if $S$ hits $\calL_n$, a
cone-point is created. Therefore, let
$\calC_\delta$ denote the event in which the trajectory of $S$ hits
$\calL_n$ at least $\delta n$ times after time $n=0$. Using {(\ref
{eqrepresRWbiss}) and keeping only configurations with empty boundary
clusters, $\gamma^{\mathsf b}=\gamma^{\mathsf f}=\varnothing$,}
\[
e^{\xi_pn}\bbP_{p}(\calM_\delta)
\geq
\rwP_{0}(\calC_\delta\cap\calR_n),
\]
where $\calR_n:=\calR(n\unitvect)$.
Dividing by $e^{\xi_pn}\bbP_{p}^{(n)}(0\conn n\unitvect)\leq e^{\xi
_pn}\bbP_{p}(0\conn n\unitvect)$ and using (\ref{eqOZsimple}) and
(\ref{eqOZasympt}), we get
%
\begin{equation}
\bbP_p^{(n)}(\calM_\delta|0\conn n\unitvect)\geq
c_{25}\rwP_0(\calC_\delta|\calR_n) ,
\end{equation}
where $c_{25}>0$ does not depend on $n$.
The next step is to express $\rwP_0(\calC_\delta|\calR_n)$ in terms of
$S^\parallel$ and $S^\bot$.
Let $\tau_0:=0$, and for $k\geq1$,
$\tau_{k}:=\inf\{m>\tau_{k-1}:S_m\in\calL\}$.
Using
(\ref{eqOZasympt}) we infer that for all $n$ and $k\leq n/2$,
$\rwP_0(\calR_{n-k})/\rwP_0(\calR_n)\geq c_{26}$ for some $c_{26}>0$,
and so by the strong Markov property,
\begin{eqnarray*}
\rwP_0(\calC_\delta|\calR_n)&=&
\rwP_0\bigl(S^\|_{\tau_{\lfloor\delta n\rfloor}}\leq n|\calR_n\bigr)\\
&=& \sum_{k\leq n/2}\rwP_0\bigl(S^\|_{\tau_{\lfloor\delta n\rfloor
}}=k\bigr)\rwP
_0(\calR_{n-k})/\rwP_0(\calR_n)\\
&\geq& c_{26}\rwP_0\bigl(S^\|_{\tau_{\lfloor\delta n\rfloor}} \leq n/2 \bigr) .
\end{eqnarray*}
Let $\bar n:=n\rwE[Y_1^\parallel]/4$. If $\calN^\parallel_n:=\max\{
k\leq n\dvtx S_k^\parallel\leq n/2\}$
denotes the number of steps performed by $S$ before leaving the strip
$\calS_{0,n\unitvect/2}$,
\begin{eqnarray*}
\rwP_0\bigl( S^\|_{\tau_{\lfloor\delta n\rfloor}} \leq n/2\bigr)&\geq&\rwP_0\bigl(
S^\|_{\tau_{\lfloor\delta n\rfloor}} \leq n/2 ;
\calN^\parallel_n\geq\bar n\bigr)
\\
&\geq& \rwP_0\bigl(L^\perp(\bar n)\geq\delta n,\calN^\parallel_n\geq\bar n\bigr)
\end{eqnarray*}
with $L^\perp(\bar n) = \#\{0\leq i\leq\bar n \dvtx S^\perp_i=0 \}$.
By an elementary large deviation estimate, $\rwP_0(\calN^\parallel_n<
\bar n)\leq e^{-c_{27}n}$ for some $c_{27}>0$.
Therefore,
\begin{eqnarray*}
\rwP_0\bigl(L^\perp(\bar n)\geq\delta n,\calN^\parallel_n\geq\bar n\bigr)
&\geq& \rwP_0\bigl(L^\perp(\bar n)\geq\delta n\bigr)-e^{-c_{27}n}\\
&=& \rwP_0\bigl(L^\perp(\bar n)\geq\delta_* \bar n\bigr)-e^{-c_{27}n} ,
\end{eqnarray*}
where $\delta_*=4\delta/\rwE[Y_1^\parallel]$.
The event $\{L^\perp(\bar n)\geq\delta_* \bar n\}$ depends only on the
transverse component $S^\perp$, which lies in $\bbZ^{d-1}$.
It follows from Corollary~\ref{corRWpinning} in Appendix
\ref{SecPinningRw} that
\[
\rwP_0 \bigl(L^\perp(\bar n) \geq\delta_* \bar n\bigr) \geq
\cases{
e^{-c \delta_*^2 n} &\quad$(d=2)$,\vspace*{2pt}\cr
e^{-c (\delta_*/|{\log\delta_*}|) n} &\quad$(d=3)$.}
\]
This proves Proposition~\ref{PropEntropyCost}.

\section{\texorpdfstring{Proof of Theorem \protect\ref{thmmain2}}{Proof of Theorem 1.4}}\label{SecPureExpon}
In this section we prove Theorem~\ref{thmmain2}: when $p'>p_c'$,
$\bbP
_{p,p'}(0\conn n\unitvect)$ has a purely exponential decay. The
underlying mechanism is
that when $\xi_p>\xi_{p,p'}$, a typical
cluster $C_{0, n\unitvect}$ connecting $0$ to $n\unitvect$ is pinned on
$\calL_n$, in the sense that it has a number of cone-points \textit{on
$\calL_n$} that grows linearly with $n$.
Cone-points of $C_{0, n\unitvect}$ lying on $\calL_n$ will be called
cone-renewals.
%
\begin{theorem}\label{propconeRenewal}
If $p'>p_c'$, then there exist $\delta= \delta(p, p')>0 $ and
$\nu_4 = \nu_4(p, p' )>0$ such that for any large enough $n$,
%
\begin{equation}
\label{eqconeRenewal}
\bbP_{p, p'}( C_{0, n\unitvect}\mbox{ has less than $\delta n$
cone-renewals}\,|
0\conn n\unitvect) \leq e^{-\nu_4 n} .\vadjust{\goodbreak}
\end{equation}
\end{theorem}

With this piece of information, irreducible components with both
endpoints on $\calL$
can be defined,
and a fairly standard renewal argument leads to the pure exponential decay.
(Note, however,
that at this point we do not even know whether under $\bbP_{p, p'}$ the cluster
$C_{0, n\unitvect}$ contains cone-points at all.)\looseness=-1

The presence of cone-points on $\calL$ will also allow to complete the
proof of Theorem~\ref{thmmain}: we show in Section
\ref{Ssecstrictmon} that
$p'\mapsto\xi_{p,p'}$ is strictly decreasing on $(p_c',1)$, and in
Section~\ref{Ssecstrictmon} that it is real analytic on the same
interval.

Assume $p'>p_c'$, and let
\[
\tau:= \xi_p-\xi_{p,p'}>0 .
\]
To prove Theorem~\ref{propconeRenewal}, we will first show that
$C_{0,n\unitvect}$ typically stays in a vicinity of size $|{\log\tau}
|/\tau$ of $\calL_n$. This implies, by a finite-energy argument, that
$C_{0,n\unitvect}$
is made of many stretches on which cone-renewals occur with positive
probability.

\subsection{Excursions away from $\calL$}
To any realization of $\{0\conn n\unitvect\}$,
we associate the smallest self-avoiding path $\pi\dvtx0\to n\unitvect$
contained in $C_{0,n\unitvect}$, as in Section~\ref{Subsecpaths}:
$\pi=(\pi_0,\pi_1,\ldots,\pi_{|\pi|})$, with $\pi_0=0$, $\pi
_{|\pi
|}=n\unitvect$.

\begin{figure}[b]

\includegraphics{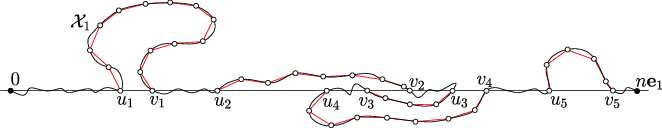}

\caption{Coarse-graining the excursions of a path $\pi\dvtx0\to
n\unitvect$.}
\label{Figexcursions}
\end{figure}

Let $K\geq1$, which will be chosen later as a function of $\tau$. Let also
\[
\tau_K(s) := \inf\{t>s \dvtx\pi_t\notin B_K(\pi_s) \} .
\]
We associate to $\pi$ a set of disjoint
pairs $(u_1,v_1),\ldots,(u_m,v_m)$ of points lying on~$\calL$, as
follows; see Figure~\ref{Figexcursions}. Let $t_0:=0$, and set, for
$j\geq1$,
\begin{eqnarray*}
s_j &:=& \inf\{s\geq t_{j-1} \dvtx\pi(s)\in\calL, \pi[s+1,\tau
_K(s)] \cap\calL= \varnothing\} ,\\
t_j &:=& \inf\{t>s_j \dvtx\pi_t\in\calL\} .
\end{eqnarray*}
We call the subpath
$\calX_j:=\pi[s_j,t_j]$ an excursion, starting at $u_j:=\pi_{s_j}$ and
ending at $v_j:=\pi_{t_j}$.

We further coarse-grain each excursion $\calX_j$ on the scale $K$. Let
$u_j^0:=u_j$ and, for $k\geq0$,
\[
u_j^{k+1} := \pi_{\tau_K(u_j^k)}.
\]
If $m_j:=\max\{k\dvtx u_j^k\in\calX_j\}$, we
call $\#_K\calX_j:= m_j$ the length of the excursion~$\calX_j$
(measured by the number of increments of size $K$).
The set of points $(u_j^0\equiv u_j,u_j^1,\ldots,u_j^{m_j},v_j)$ is
called the skeleton of $\calX_j$.
Sometimes, $u_j^{{m_j}}\equiv v_j$, but in all cases,
$|u_j^{m_j}-v_j|\leq K$.

We denote by $\{u\stackrel{m}{\curvearrowright} v\}$ the event in which
there exists a path which is an excursion of length $m$ starting at $u$
and ending at $v$.
%
\begin{lem}\label{lemexcursion}
There exists $K_0=K_0(\tau)$ and $c_{28}=c_{28}(\tau)>0$ such that
if $K\geq K_0$,
\[
\bbP_{p,p'}\bigl(u\stackrel{m}{\curvearrowright} v\bigr) \leq
e^{-\xi
_{p,p'}|v-u| - c_{28}\tau K m}.
\]
\end{lem}
\begin{pf}
Denote by $\calX$ any excursion occurring in $\{u\stackrel
{m}{\curvearrowright} v\}$. That is, $\#_K\calX=m$. Let $(u^0,\ldots
,u^m)$ be a skeleton, where for the sake of simplicity, we assume that
$u^m=v$. By construction,
the event
\[
\{ \calX\mbox{ has skeleton $(u^0,\ldots,u^m)$}\}
\]
implies that
there are disjoint connections $u^0 \stackrel{\calL^c}{\conn} u^1,
\ldots, u^{m-1} \stackrel{\calL^c}{\conn} u^m$.
By the BK inequality,
\begin{eqnarray*}
\bbP_{p,p'} \bigl(\calX\mbox{ has skeleton $(u^0,\ldots,u^m)$}\bigr)
&\leq&
\prod_{i=1}^{m} \bbP_{p,p'}\bigl(u^{i-1}\stackrel{\calL^c}{\conn}
u^{i}\bigr)\\
&\leq&
\prod_{i=1}^{m} \bbP_{p}(u^{i-1}\conn u^{i})\\
&\leq&\prod_{i=1}^{m}e^{-\xi_p(u^i-u^{i-1})} .
\end{eqnarray*}
If $z\in\bbR^d$, let $k\in\{1,\ldots,d\}$ be such that
$\langle\mathbf{e}_k,z\rangle={\max_{k'}}|\langle\mathbf
{e}_{k'},z\rangle|$. Then, using (\ref{eqWulff}) and since $\xi
_p\mathbf{e}_k\in\partial W_p$,
\begin{eqnarray*}
\xi_p(z)&=&\sup_{t\in\partial W_p} \langle t,z\rangle
\geq\xi_p \langle\mathbf{e}_k,z\rangle
= \xi_{p,p'}\langle\mathbf{e}_k,z\rangle+\tau\langle\mathbf
{e}_k,z\rangle\\
&\geq& \xi_{p,p'} |\langle\unitvect,z\rangle| + c_{29}\tau|z|
\end{eqnarray*}
for some constant $c_{29}=c_{29}(d)>0$. Since
\[
\sum_{i=1}^m |\langle\unitvect,u^i-u^{i-1}\rangle| \geq|v-u|
\]
and $|u^i-u^{i-1}|\geq K$ for all $i=1,\ldots,m$, we get
\[
\bbP_{p,p'}\bigl(\calX\mbox{ has skeleton $(u^i)_{i=0,\ldots
,m}$}\bigr)
\leq e^{-\xi_{p,p'}|v-u|- c_{29}\tau K m}.
\]
When $u^m\neq v$, a similar computation leads to the same bound.
Since the number of skeletons with $m$ increments is $O
((K^{d-1})^mK^d)$, the conclusion follows by taking $K\geq K_0$, with
$K_0$ large enough in order that $\frac{\log K_0}{K_0}$ be sufficiently
small compared to $\tau$.\vadjust{\goodbreak}
\end{pf}

Let $\#_K\pi:= \sum_{j}\#_K\calX_j$ denote the total number of
increments in the excursions of $\pi$.
%
\begin{prop}\label{Proexcpi} Let $0<\varepsilon<1$.
There exists $K_1=K_1(\tau,\varepsilon)$ and $c_{30}=c_{30}(\varepsilon)>0$
such that for all $K\geq K_1$,
%
\begin{equation}\label{equppertotalexc}\bbP_{p,p'} ( \#_K\pi
\geq
\varepsilon n/K | 0\conn n\unitvect)
\leq e^{-c_{30} n}.
\end{equation}
\end{prop}
\begin{pf}
For a {collection of triples} $(u_j,v_j,m_j)_{j=1}^M$, let $\calP
((u_j,v_j,m_j)_{j=1}^M)$ denote the event on which there exists a path
$\pi\dvtx0\to n\unitvect$ with $M$ excursions, the $j$th excursion $\calX
_j$, starting at $u_j$ and ending at $v_j$, and being such that $\#
_K\calX_j=m_j$.
The event $\calP((u_j,v_j,m_j)_{j=1}^M)$ implies the disjoint occurrence
\[
\{ v_{0}\conn u_1\}\circ
\bigl\{u_1\stackrel{m_1}{\curvearrowright} v_1\bigr\}
\circ\cdots\circ
\{ v_{M-1}\conn u_M\}\circ
\bigl\{u_M\stackrel{m_M}{\curvearrowright} v_M\bigr\}.
\]
Assuming $K$ is larger than the $K_0$ of Lemma~\ref{lemexcursion}, and
by the BK inequality,
\begin{eqnarray*}
\bbP_{p,p'}(\calP((u_j,v_j,m_j)_{j=1}^M))
&\leq&
\prod_{j=1}^{M}\bbP_{p,p'}(v_{j-1}\conn u_j)\bbP_{p,p'}\bigl(u_j\stackrel
{m_j}{\curvearrowright} v_j\bigr)\\
&\leq&
\prod_{j=1}^{M} e^{-\xi_{p,p'} (|u_j-v_{j-1}|+|v_j-u_j|)}
e^{-c_{28}\tau K m_j} .
\end{eqnarray*}
We then sum over the triples $(u_j,v_j,m_j)_{j=1}^M$.
Denote by $I\supset\calL_n$ the smallest interval of $\calL$
containing all the points $u_j,v_j$, $j=1,\ldots,M$.
Observe that $\sum_{j=1}^{M}(|u_j-v_{j-1}|+|v_j-u_j|)\geq|I|$.
We first sum over the possible positions of $I$, then over the
positions of the $M\geq1$
distinct points $u_j$ in $I$, then over the $m_j$'s satisfying $\sum
_{j=1}^Mm_j\geq\varepsilon n/K$,
and finally over the endpoints $v_j$.
Since to a given point $u_j$ correspond at most $2K(m_j+1)$ points $v_j$,
\begin{eqnarray*}
&&\bbP_{p,p'} ( \#_K\pi\geq\varepsilon n/K, 0\conn n\unitvect)\\
&&\qquad\leq
\sum_{I\supset\calL_n} e^{-\xi_{p,p'} |I|} \sum_{M\geq1} \pmatrix
{|I|\cr M}
\mathop{\sum_{m_1,\ldots,m_M \geq1}}_{\sum_{j} m_j\geq\varepsilon n/K}
\prod_{j=1}^M \bigl(2K(m_j+1)\bigr)
e^{- c_{28}\tau K m_j}.
\end{eqnarray*}
We choose $K_1\geq K_0$ large enough so that, for all $K\geq K_1$ and
all $m\geq1$,
$(2K(m+1))e^{- c_{28}\tau K m} \leq e^{- c_{31}\tau K m}$.
Proceeding as on page \pageref{pagenummethod},
%
\begin{equation}\label{eqboundnumbincr}
\sum_{M\geq1} \pmatrix{|I|\cr M} \mathop{\sum_{m_1,\ldots,m_M
\geq
1}}_{\sum_{j} m_j\geq\varepsilon n/K} \prod_{j=1}^M
e^{- c_{31}\tau K m_j}\leq e^{c_{32}e^{-c_{31}\tau K/2}|I|} .
\end{equation}
Then, notice that there are $\ell-n$ intervals $I\supset\calL_n$ of
fixed length $|I|=\ell$.
Therefore, summing over $|I|$ gives
\[
\bbP_{p,p'} ( \#_K\pi\geq\varepsilon n/K, 0\conn n\unitvect)
\leq e^{-c_{33}\tau\varepsilon n}e^{-\xi_{p,p'}n} .\vadjust{\goodbreak}
\]
Since $\bbP_{p,p'}(0\conn n\unitvect)=e^{-\xi_{p,p'}(1+o(1))n}$ as
$n\to\infty$,
we get (\ref{equppertotalexc}) once $K$ is sufficiently
large.
\end{pf}

We then turn to the study of the deviations of $C_{0,n\unitvect}$ away
from its smallest connecting path $\pi\subset C_{0,n\unitvect}$.

Let $\pi$ be a given path, which we here consider together with its set
of edges.
Let
$R_0:=\max\{|z-0|\dvtx0\stackrel{\pi^c}{\conn}z\}$ and $z_0\in
C_{0,n\unitvect}$ be any point at which the max is attained.
Let $\hat\pi^0$ be the smallest path realizing the connection between
$0$ and $z_0$, disjoint from $\pi$.
Inductively, for $s= 1,\ldots,|\pi|$,
let $\Pi_s:=\bigcup_{0\leq t<s}\hat\pi^t$,
\[
R_s:=\max\bigl\{|z-\pi_s|\dvtx\pi_{s}
\stackrel{{(\pi\cup\Pi_s)^c}}{\longleftrightarrow}z\bigr\},
\]
$z_s\in C_{0,n\unitvect}$ be any point at which the max is attained, and
$\hat\pi^s$ be any path realizing the connection between $\pi_{s}$ and
$z_s$, disjoint from {$\pi\cup\Pi_s$}.
%
\begin{prop}\label{Prohairs}
Let $0<\varepsilon<1$. There exists $K_3=K_3(p,p',\varepsilon)$ and
$c_{34}=c_{34}(p,p',\varepsilon)>0$ such that if $K\geq K_3$, then
%
\begin{equation}\label{eqboundhairs}\bbP_{p,p'}\Biggl(\sum
_{s=0}^{|\pi
|}R_s\mathbf{1}_{\{{R_s\geq K}\}}\geq\varepsilon n\Big|0\conn
n\unitvect\Biggr)
\leq e^{-c_{34}n} .
\end{equation}
\end{prop}
\begin{pf}
We know from Proposition~\ref{Proexcpi} that under $\bbP
_{p,p'}(\cdot
|0\conn n\unitvect)$,
the number of increments of the skeleton of a typical path $\pi\dvtx0\to
n\unitvect$ is at most $\varepsilon n/K$. We can therefore assume, in
particular, that
%
\begin{equation}\label{eqboundsizepi}|\pi|\leq c_{35}K^{d-1}n .
\end{equation}
For a fixed path $\pi$, let $\calF_\pi$ denote the event in which
$\pi$
is the smallest self-avoiding path connecting $0$ to $n\unitvect$.
Arguing as for (\ref{eqcoupling}), we get
$\bbP_{p,p'}(\cdot| \calF_\pi)
\preccurlyeq\bbP_{p,p'}(\cdot)$ on~$\pi^c$.
The event $\{ \sum_{s=0}^{|\pi|}R_s\mathbf{1}_{\{{R_s\geq K}\}
}\geq\varepsilon
n\}$ is
$\pi^c$-measurable and increasing.
Therefore, by the BK inequality and Lemma
\ref{lemexpdecay},
\begin{eqnarray*}
&&\bbP_{p,p'}\Biggl(\sum_{s=0}^{|\pi|}R_s\mathbf{1}_{\{{R_s\geq K}\}
}\geq\varepsilon
n\Big
|\calF_\pi\Biggr)\\
&&\qquad
\leq\bbP_{p,p'}\Biggl(\sum_{s=0}^{|\pi|} R_s\mathbf{1}_{\{{R_s\geq
K}\}}\geq
\varepsilon
n\Biggr)
\\
&&\qquad\leq\mathop{\mathop{\sum_{r_1,\ldots,r_{|\pi|}:}}_{\sum_sr_s\geq
\varepsilon n,}}_{r_s\geq K}
\mathop{\sum_{z_1,\ldots,z_{|\pi|}:}}_{|z_s-\pi_s|=r_s}
\bbP_{p,p'}\bigl(
\{\pi_0\conn z_0\}\circ\cdots\circ\bigl\{\pi_{|\pi|}\conn z_{|\pi|}\bigr\}
\bigr)\\
&&\qquad\leq\mathop{\mathop{\sum_{r_1,\ldots,r_{|\pi|}:}}_{\sum_sr_s\geq
\varepsilon n,}}_{r_s\geq K}
\prod_{s=0}^{|\pi|}c_{36}r_s^{d-1}e^{-c_{11}\xi_p^*r_s} .\vadjust{\goodbreak}
\end{eqnarray*}
The proof then follows the same lines as before: if $K$ is large
enough, then
$c_{36}r^{d-1}e^{-c_{11}\xi_p^*r}\leq e^{-c_{37}r}$ for all $r\geq K$.
The summation can thus be done as in (\ref{eqboundnumbincr}), and
using (\ref{eqboundsizepi}) gives (\ref{eqboundhairs}).
\end{pf}

Let $\calT_{2K}$ be the tube containing points whose {Euclidean}
distance to $\calL$ is $\leq2K$, and
consider the cone
\[
\calY:=\{x\dvtx\langle x,\unitvect\rangle\geq|x^\perp|\} .
\]
For each $x\in\calL^c$, let $z_+$ (resp., $z_-$) denote the largest
(resp., smallest) point of $\calL$ such that $x\in z_+-\calY$ (resp.,
$x\in z_-+\calY$).
The segment $[z_-,z_+]$ is called the shade of $x$.
Let $\calS_n\subset\calL_n$ be the set of points of $\calL_n$ who lie
in the shade of at least one point of $C_{0,n\unitvect}\cap(\calT_{2K})^c$.
The points of $\calR_n:=\calL_n\setminus\calS_n$ are candidates for
being cone-renewals.

It is easy to see that
\[
|\calS_n|\leq c_{38} K\sum_j\#_K\calX_j+c_{39}K\sum_{s=0}^{|\pi
|}R_s\mathbf{1}_{\{{R_s\geq K}\}} ,
\]
where $c_{38}$ and $c_{39}$ depend only on the dimension $d$.
As a corollary of Propositions~\ref{Proexcpi} and~\ref{Prohairs},
$|\calS_n|=O(\varepsilon n)$. More precisely, for a fixed $0<\eta<1$, $K$
can be taken large enough so that
%
\begin{equation}\label{eqprerenewals}
\bbP_{p,p'}\bigl(|\calR_n|\geq(1-\eta) n|0\conn n\unitvect
\bigr)\geq
1-e^{-c_{40}n}
\end{equation}
with $c_{40}>0$ depending on $p$, $p'$ and $\eta$.
\begin{pf*}{Proof of Theorem~\ref{propconeRenewal}}
We apply a local surgery under
$\bbP_{p,p'}(\cdot|0\conn n\unitvect)$, to show that $\calL_n$
contains many cone-renewals (see Figure~\ref{figblockNew}).
Consider the partition of $\calT_n$ into neighboring disjoint blocks
%
\begin{figure}

\includegraphics{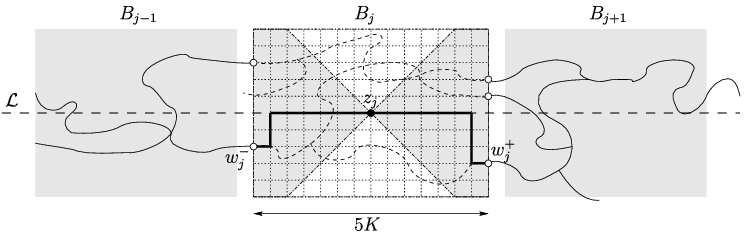}

\caption{The local surgery inside the block $B_j$, turning a
pre-renewal $z_j$ into a cone-renewal: open a minimal path (in bold)
connecting $w_j^-$ to $w_j^+$, and close all other edges of~$B_j$ (dotted).
By the finite energy property, this event has probability $p_*=e^{-O(K^d)}>0$.}
\label{figblockNew}\vspace*{6pt}
\end{figure}
$B_j$ of lengths $5K$, centered at points~$z_j$. If $z_j\in\calR_n$, we
call $z_j$ a pre-renewal. Assume $z_j$ is a pre-renewal. Let $F_j^{-}$,
$F_j^+$ denote the two faces of $B_j$ which are orthogonal to
$\calL_n$, and let $W_j^-\subset F_j^-$ (resp., $W_j^+\subset F_j^+$)
denote the points of $C_{0,n\unitvect}\cap F_j^-$ (resp.,
$C_{0,n\unitvect}\cap F_j^+$) which are connected to $0$ (resp.,
$n\unitvect$) by a path not intersecting $B_j$. Let $w_j^\pm$ denote
the smallest point (in lexicographical order) of $F_j^\pm$. Under
$\bbP_{p,p'}(\cdot|w_j^-,w_j^+,0\conn n\unitvect)$, independently of
the edges living outside $B_j$, $w_j^-$ is connected to $w_j^+$ by a
minimal path going through $z_j$, turning $z_j$ into a cone-renewal
with positive probability, bounded below by some $p_*>0$ depending on
$K$.

The variables $X_i:=\mathbf{1}_{\{{z_i\ \mathrm{is}\ \mathrm{a}\
\mathrm{cone}\mbox{-}\mathrm{renewal}}\}}$
can thus be
coupled to i.i.d. Ber\-noulli
variables $Y_i$ of parameter $p$, giving
\[
\bbP_{p,p'}\Biggl(\sum_{i=1}^{n/5K}X_i\leq p_*n/(10 K)
\Big|0\conn n\unitvect\Biggr)\leq
\rwP\Biggl(\sum_{i=1}^{n/5K}Y_i\leq p_*n\Big/(10K)\Biggr)\leq e^{-c_{41}n} .
\]
Together with (\ref{eqprerenewals}), this proves the claim.
\end{pf*}

Let us complete the proof of Theorem~\ref{thmmain2}.
We first define the irreducible components $\zeta_j$ of
$C_{0,n\unitvect
}$, which are cone-confined and which, in contrast to the $\gamma_j$ of
Section~\ref{SecRWREP}, have both their endpoints on $\calL_n$.

Let us denote by $\{w_1,\ldots,w_{m+1}\}\subset C_{0,n\unitvect}$ the
cone-renewals that lie on $\calL_n$, ordered according to their first
component.
By Theorem~\ref{propconeRenewal}, $m$ is typically of order $n$.
The subgraphs
\[
\zeta_j:= C_{0,n\unitvect}\cap\calS_{w_{j}, w_{j+1}}
\]
are called cone-confined irreducible components of $C_{0,n\unitvect}$.
The complement
{$C_{0,n\unitvect}\setminus(\zeta_1\cup\cdots\cup\zeta_m)$} can
contain, at most, two connected components.
If it exists, the component containing $0$ (resp., $n\unitvect$) is
denoted $\zeta^{\mathsf b}$ (resp., $\zeta^{\mathsf f}$), and called backward
(forward) irreducible.
Keeping in mind that we are here working with the cone $\calY$ rather
than $\calY^>$ and that the edges on $\calL$ are opened with
probability~$p'$, all
the definitions of Section~\ref{SecRWREP} extend with almost no
changes to
the irreducible components $\zeta$.
In particular, we can define independent events $\Xi^{\mathsf b},\Xi
_1,\ldots
,\Xi_m,\Xi^{\mathsf f}$ so that
\[
\bbP_{p,p'}(C_{0,n\unitvect}=\zeta^{\mathsf b}\sqcup\zeta_1\sqcup
\cdots\sqcup
\zeta_m\sqcup\zeta^{\mathsf f})=\bbP_{p,p'}(\Xi^{\mathsf b})\Biggl(\prod
_{j=1}^m\bbP_{p,p'}(\Xi_j)\Biggr)\bbP_{p,p'}(\Xi^{\mathsf f}) .
\]
One can thus define, for $u\geq1$, $v\leq-1$,
\[
\rho_{\mathsf b}'(u):=e^{\xi_{p,p'}u}\mathop{\sum_{\zeta^{\mathsf b}\ni
0:}}_{{\mathsf b}(\zeta^{\mathsf b})=u}\bbP_{p,p'}(\Xi^{\mathsf b})
,\qquad
\rho_{\mathsf f}'(v):=e^{\xi_{p,p'}|v|}\mathop{\sum_{\zeta^{\mathsf f}\ni
0:}}_{{\mathsf f}(\zeta^{\mathsf f})=v}\bbP_{p,p'}(\Xi^{\mathsf f}) .
\]
{By (\ref{eqconeRenewal})}, these weights satisfy the following bounds:
%
\begin{equation}\label{eqproprnewweights}
\rho_{\mathsf b}'(u)\leq e^{-\nu_4|u|} ,\qquad
\rho_{\mathsf f}'(v)\leq e^{-\nu_4|v|} .
\end{equation}
Moreover, $q'(\ell):=e^{\xi_{p,p'}\ell}f_\ell$ with
%
\begin{equation}\label{eqdeffl}
f_\ell:=
\mathop{\sum_{\zeta_1\ni0:}}_{{\mathsf f}(\zeta_1)=0,{\mathsf b}(\zeta
_1)=\ell
}\bbP_{p,p'}(\Xi_1)
\end{equation}
defines a probability distribution on $\bbN$.
{Again, by (\ref{eqconeRenewal}),}
%
\begin{equation}\label{eqboundfl}
f_\ell\leq e^{-\xi_{p,p'}\ell-\nu_4\ell} ,
\end{equation}
which implies
%
\begin{equation}\label{eqborneqprime}
q'(\ell)\leq e^{-\nu_4\ell} .
\end{equation}
Up to a term of order $e^{-\nu_4 n}$ [compare with (\ref
{eqidentwithincrements})],
%
\begin{equation}\label{eqdecompRenewal}\quad
e^{\xi_{p,p'}n}\bbP_{p,p'}(0\conn n\unitvect)=\sum_{u,v}\rho_{\mathbf
b}'(u)\rho_{\mathbf f}'(v)\sum_{m\geq1}
\mathop{\sum_{\ell_1,\ldots, \ell_{m}:}}_{\sum_j\ell_j=n+v-u}
\prod_{j=1}^m q'(\ell_j)
.
\end{equation}
As before, due to (\ref{eqproprnewweights}), the sum in (\ref
{eqdecompRenewal})
can be restricted to those $u,v$ that satisfy $|u|,|v|\leq
n^{1/2-\alpha
}$, for some small $\alpha>0$.
Let thus $\tau_k$, $k\geq1$, be an i.i.d. sequence with distribution
$Q'(\tau_1=\ell):=q'(\ell)$. Then, (\ref{eqdecompRenewal}) writes
\begin{eqnarray*}
&&e^{\xi_{p,p'}}
\bbP_{p,p'}(0\conn n\unitvect)\\
&&\qquad=\sum_{u,v}
\rho_{\mathbf b}'(u)\rho_{\mathbf f}'(v)Q'\Biggl(\exists m\geq1\mbox{ such that
}\sum_{j=1}^m\tau_j=n+v-u\Biggr) .
\end{eqnarray*}
By (\ref{eqborneqprime}), $E_{Q'}[\tau_1]<\infty$. Moreover,
$q'(\ell
)>0$ for all $\ell\geq1$, and therefore, by the renewal theorem,
\[
Q'\Biggl(\exists m\geq1\mbox{ such that }\sum_{j=1}^m\tau
_j=n+v-u
\Biggr)\to\frac{1}{E_{Q'}[\tau_1]}
\]
as $n\to\infty$, uniformly in $|u|,|v|\leq n^{1/2-\alpha}$.
This proves Theorem~\ref{thmmain2}.

\subsection{\texorpdfstring{Strict monotonicity of $p'\mapsto\xi_{p,p'}$}
{Strict monotonicity of p' to xi p,p'}}
\label{Ssecstrictmon}

Assume $p'>p'_c$, that is, $\xi_{p,p'}<\xi_p$.
Consider the measures $\bbP_{p,p'}^{(n)}$ defined in Section~\ref{secLowerB}.
If $a_n\gg n$ is taken large enough, then
we can write $\xi_{p,p'}=\lim_{n\to\infty}\xi_{p,p'}^{(n)}$, where
\[
\xi_{p,p'}^{(n)}:=-\frac{1}{n}\log\bbP_{p,p'}^{(n)}(0\conn
n\unitvect
) .
\]
Therefore,
\[
\frac{\dd\xi_{p,p'}^{(n)}}{\dd p'}=-\frac{1}{n}\frac{\dd\bbP_{p,p'}^{(n)}(0\conn n\unitvect)/{\dd
p'}}{\bbP_{p,p'}^{(n)}(0\conn
n\unitvect)}.
\]
By Theorem~\ref{propconeRenewal}, the expected number of cone-renewals
under $\bbP_{p,p'}(\cdot| 0\conn n\unitvect)$ grows linearly with $n$.
Since each cone-renewal is adjacent to two edges which are pivotal for
$\{0\conn n\unitvect\}$, we can use Russo's Formula as before to find a
constant $c_{42}>0$ such that
\[
\frac{1}{n}
\frac{\dd\bbP_{p,p'}^{(n)}(0\conn n\unitvect)/{\dd p'}}{\bbP
_{p,p'}^{(n)}(0\conn n\unitvect)}\geq
c_{42}.
\]
This implies that $\frac{\dd\xi_{p,p'}^{(n)}}{\dd p'}\leq-c_{42}$,
uniformly in $n$. $p'\mapsto\xi_{p,p'}$ is therefore strictly
decreasing on $(p_c',1)$, since for all $p_c'<p_1'<p_2'<1$,
\[
\xi_{p,p_2'}-\xi_{p,p_1'}=\lim_{n\to\infty}\bigl(\xi
_{p,p_2'}^{(n)}-\xi
_{p,p_1'}^{(n)}\bigr)
=\lim_{n\to\infty}\int_{p_1'}^{p_2'}\frac{\dd\xi
_{p,p'}^{(n)}}{\dd
p'}\,\dd p'\leq
-c_{42}(p_2'-p_1') .
\]

\subsection{\texorpdfstring{Analyticity of $p'\mapsto\xi_{p,p'}$}
{Analyticity of p' to xi p,p'}}\label{Ssecanalyticity}

Fix $p<p_c$. Consider $f_\ell=f_\ell(p')$ defined in~(\ref{eqdeffl}).
Observe that $f_\ell$ can be put in the form of a polynomial in $p'$,
$f_\ell(p')=\sum_{k=0}^\ell a_k^{(\ell)}p'^k$, with $a_k^{(\ell
)}\geq0$. It can therefore be continued as an analytic function
$w\mapsto f_\ell (w)$ in the complex plane. Let
\[
\Phi(w,z):=\sum_{\ell\geq1}f_\ell(w)e^{z\ell} .
\]
Since
\[
\Phi(p',\xi_{p,p'})=\sum_{\ell\geq1}f_\ell(p')e^{\xi_{p,p'}\ell
}=\sum
_{\ell\geq1} q'(\ell)=1 ,
\]
the analyticity of $p'\mapsto\xi_{p,p'}$ will follow by solving $\Phi
(w,z)=1$ for $z$, in a neighborhood of $(p',\xi_{p,p'})$.
To do so, we must verify that $(w,z)\mapsto\Phi(w,z)$ is analytic in a
domain of $\bbC^2$ containing $(p',\xi_{p,p'})$, and that
$\frac{\partial\Phi}{\partial z}|_{(p',\xi_{p,p'})}\neq0$.
If $w\in D_\delta(p'):=\{w\in\bbC\dvtx|w-p'|<\delta\}$,
\[
|f_\ell(w)|\leq\sum_{k=0}^\ell a_k^{(\ell)}|w|^k\leq\sum
_{k=0}^\ell
a_k^{(\ell)}(p'+\delta)^k
\leq(1+\delta/p')^\ell f_\ell(p') .
\]
We can therefore choose $\delta=\delta(p,p')>0$ small enough to
ensure that
\[
\sup_{w\in D_\delta(p')}\biggl|\frac{f_\ell(w)}{f_\ell(p')}\biggr|
\leq e^{\nu_4 \ell/3} .
\]
We also take $\varepsilon>0$ such that $\sup_{z\in D_\varepsilon(\xi
_{p,p'})}|e^{z\ell}|\leq e^{(\xi_{p,p'}+\nu_4 /3)\ell}$.
Remembering the bound for $f_\ell(p')$ in (\ref{eqboundfl}), we thus get
\[
\sum_{\ell\geq1}\sup_{w\in D_\delta(p')}\sup_{z\in D_\varepsilon
(\xi
_{p,p'})}|f_\ell(w)e^{z\ell}|<\infty.
\]
Therefore, $\Phi$ defines an analytic function of $(w,z)$ in the polydisc
$D_\delta(p')\times D_\varepsilon(\xi_{p,p'})$. Moreover,
\[
\frac{\partial\Phi}{\partial z}\bigg|_{(p',\xi_{p,p'})}
=\sum_{\ell\geq1}\ell f_\ell(p') e^{\xi_{p,p'}\ell}>0 .
\]
The conclusion follows by the implicit function theorem.

\begin{appendix}

\section{Renewals}\label{Apprenewal}

Let $(a_n)_{n\geq0}$ and $(b_n)_{n\geq1}$ be nonnegative sequences
satisfying $a_0=1$, and the renewal equation
%
\begin{equation}\label{eqrenewalanbn}
a_n=\sum_{k=0}^{n-1}a_kb_{n-k} \qquad\mbox{for all }n\geq1 .
\end{equation}
Iterating (\ref{eqrenewalanbn}) gives
%
\begin{equation}\label{eqrenewalbasic}
a_n=\sum_{m=1}^n\mathop{\sum_{k_1,\ldots,k_m}}_{\sum_jk_j=n}\prod
_{j=1}^mb_{k_j} \qquad \mbox{for all }n\geq1 .
\end{equation}
As a consequence, in terms of the generating functions
\[
\bbA(s)=\sum_{n\geq0}a_ns^n ,\qquad \bbB(s)=\sum_{n\geq1}b_ns^n ,
\]
equation (\ref{eqrenewalanbn}) takes the form
%
\begin{equation}\label{eqrenewalGF}
\bbA(s)=\frac{1}{1-\bbB(s)} .
\end{equation}
The following classical result (or variants of it) is used at various
places in the paper.
%
\begin{lem}\label{lemrenewal}
Assume that the radii of convergence of $\bbA$ and $\bbB$, denoted,
respectively,
$r_\bbA$ and $r_\bbB$, satisfy $r_\bbB>r_\bbA>0$. Then $\bbB
(r_{\bbA})=1$.
In particular, the numbers $q_k:=b_kr_\bbA^{k}$ ($k\in\bbN$) define a\vadjust{\goodbreak}
probability distribution on $\bbN$.
Moreover, if $b_k>0$ for all $k\geq1$, then
%
\begin{equation}\label{eqrenewaltheorem}
r_\bbA^{n}a_n\to\biggl(\sum_{k\geq1}kq_k\biggr)^{-1} .
\end{equation}
\end{lem}
\begin{pf}
Since its coefficients are $\geq0$, $\bbA(s)$ is singular at $s=r_\bbA
$, and therefore (\ref{eqrenewalGF}) gives $\bbB(r_{\bbA})=1$. Let
$\tau_1,\tau_2,\ldots$ denote an i.i.d. sequence
with distribution $Q(\tau_1=k):=q_k$. Then (\ref{eqrenewalbasic})
becomes
\[
r_\bbA^{n}a_n=Q(\exists M\geq1\dvtx \tau_1+\cdots+\tau_M=n
) .
\]
By the renewal theorem,
\[
Q(\exists M\geq1\dvtx \tau_1+\cdots+\tau_M=n)\to\frac
{1}{E_Q[\tau_1]}=\frac{1}{\sum_{k\geq1}kq_k} ,
\]
which proves (\ref{eqrenewaltheorem}).
\end{pf}

\section{Pinning for a random walk}\label{SecPinningRw}
In this section, we consider the pinning problem for a random walk on
$\bbZ^{d}$.
This is a classical problem (see, e.g., the book~\cite{Giacomin} and
references therein); nevertheless, for the convenience of the reader,
we state and prove the relevant claims.
The dimension $d$ of this section corresponds to dimension
$d-1$ in the paper, since the walk $X$ introduced below is associated
to the transverse component $S^\perp$ of the random walk representation
of $C_{0,n\unitvect}$.

Consider a random walk $X=(X_n)_{n\geq0}$ on $\bbZ^d$ such that (i)
$X$ is nonlattice, (ii)~$X_0=0$, (iii) the increments $X_{i+1}-X_i$
have zero expectation and exponential tails.
We denote the law of $X$ by $\rwP$. We introduce the measure $\rwP
_N^\varepsilon$ defined by
\[
\frac{\dd\rwP_N^\varepsilon}{\dd\rwP} = \frac{e^{\varepsilon L(N)}
\mathbf{1}_{\{{X_N=0}\}}}{\rwZ_N^\varepsilon},
\]
where $L(N) = \sum_{n=1}^N \mathbf{1}_{\{{X_n=0}\}}$ is the local
time at the origin,
$\varepsilon\geq0$ is the pinning parameter, and
\[
\rwZ_N^\varepsilon= \rwE\bigl[e^{\varepsilon L(N)} \mathbf{1}_{\{{X_N=0}\}}\bigr]
\]
is the normalizing partition function.

The first result shows that in dimensions $1$ and $2$, and only in
those dimensions, an arbitrary $\varepsilon>0$ leads to an exponential
divergence of $\rwZ_N^\varepsilon$.
%
\begin{theorem}
\label{thmRWfreeEnergy} For all $d\geq1$,
there exists $\varepsilon_{c}=\varepsilon_{c}(d)\geq0$ such that
\[
f(\varepsilon) = \lim_{N\to\infty} \frac1N\log\rwZ_N^\varepsilon
\cases{
= 0, &\quad if $\varepsilon< \varepsilon_{c}$,\cr
> 0, &\quad if $\varepsilon> \varepsilon_{c}$.}
\]
In dimensions $1$ and $2$, $\varepsilon_{c}(1)=\varepsilon_{c}(2)=0$,
while $\varepsilon_{c}(d)>0$ for all $d\geq3$.
Moreover, in dimensions $1$ and $2$, there exist $c_{43},c_{44}>0$ such that
%
\begin{equation}\label{eqasymptfree}
f(\varepsilon) =
\cases{
c_{43}\varepsilon^2 \bigl(1+o(1)\bigr) &\quad$(d=1)$,\vspace*{2pt}\cr
e^{-c_{44}/\varepsilon (1+o(1))} &\quad$(d=2)$.}
\end{equation}
\end{theorem}
\begin{pf}
We omit the proof of the existence of the free energy $f(\varepsilon)$,
which is standard.
The existence of $\varepsilon_c(d)$ follows by monotonicity.
Let $\tau_0:=0$ and, for $k\geq1$,
$\tau_k:=\inf\{n>0 \dvtx X_{\tau_{k-1}+n}=0 \}$.
It is well known~\cite{Caravenna,JainPruitt} that, as $k\to\infty$,
%
\begin{equation}
\label{eqRWtauAsymp}
\rwP(\tau_1 = k) =
\cases{
c_{45} k^{-3/2} \bigl(1+o(1)\bigr) &\quad$(d=1)$,\vspace*{2pt}\cr
c_{46} k^{-1}(\log k)^{-2} \bigl(1+o(1)\bigr) &\quad$(d=2)$,}
\end{equation}
for some constant $c_{45}$ and $c_{46}=2\pi\sqrt{\det\Gamma}$, with
$\Gamma$ the covariance matrix of~$X$.
Notice now that $\rwZ_N^\varepsilon$ satisfies the following renewal equation:
\[
\rwZ_N^\varepsilon
= \sum_{k=1}^{N} e^\varepsilon\rwP(\tau_1=k) \rwZ_{N-k}^\varepsilon,
\]
where we have set $\rwZ_0^\varepsilon:=1$.
Consider the generating function
$\bbA(s):=\break\sum_{N\geq0}Z_N^\varepsilon s^N$
whose radius of convergence is given by $e^{-f(\varepsilon)}\leq1$.
Proceeding as in Appendix~\ref{Apprenewal},
%
\begin{equation} \label{eqrenF}
\bbA(s)=1/\bigl(1-\bbB(s)\bigr) ,
\end{equation}
where $\bbB(s):=\sum_{k\geq1}s^ke^\varepsilon\rwP(\tau_1=k)$.
Observe that $\bbB(s)$ converges for all $s\in[0,1]$.
Since $\bbB$ is monotone, we have $\bbB(s)\leq\bbB(1)
=e^\varepsilon\rwP(\tau_1<\infty)$ for all $s<1$.

In dimension $d\geq3$, the walk is transient: $\rwP(\tau_1<\infty)<1$.
Therefore,
if $\varepsilon< \varepsilon_{c}(d):=|{\log\rwP}(\tau_1<\infty)|$, we
have $\bbB(1)<1$, so $\bbA(s)$ converges for all $s\leq1$ and
therefore $f(\varepsilon)=0$.
Now if $\varepsilon>\varepsilon_{c}$, then $\bbB(1)=e^{\varepsilon-\varepsilon
_c}>1$. Therefore,
$\bbB(s)>1$ for $s$ sufficiently close to $1$. This implies by (\ref
{eqrenF}) that
the radius of convergence of $\bbA$ is strictly smaller than $1$, and
so $f(\varepsilon)>0$.

In dimensions $d=1,2$, the walk is recurrent: $\rwP(\tau_1<\infty)=1$.
Therefore,
$\bbB(1)=e^\varepsilon>1$
for all $\varepsilon>0$, which implies that
$\bbB(s)>1$ as soon as
$s<1$ is sufficiently close to $1$. As before, this implies that
$f(\varepsilon)>0$.
Therefore, $\varepsilon_c(1)=\varepsilon_c(2)=0$.
Since $f(\varepsilon)$ is characterized by the unique number $f>0$ for which
$\bbB(e^{-f})= 1$, that is,
\[
\sum_{k\geq1}e^\varepsilon\rwP(\tau_1=k)e^{-f k}= 1.
\]
Using (\ref{eqRWtauAsymp}), an integration by parts in this last sum
shows that
as $\varepsilon\downarrow0$, $f(\varepsilon)$ behaves as in (\ref{eqasymptfree}).
\end{pf}

The second theorem provides some information about the local time at
the origin under $\rwP_N^\varepsilon$.
%
\begin{theorem}\label{thmRWpinning}
Assume that $d=1$ or $2$, and $\varepsilon>0$.
Let $\hat\tau_k$ be an i.i.d. sequence with distribution $Q(\hat
\tau
_1=k):=e^\varepsilon\rwP(\tau_1=k)e^{-f(\varepsilon)}$.
Then for all \mbox{$\eta>0$},
%
\begin{equation}
\label{eqLLNRW}
\rwP_N^\varepsilon\biggl( \biggl|\frac{L(N)}N - \frac1{\trwE[\hat
\tau_1]}\biggr|\geq\eta\biggr) \to0 .
\end{equation}
Moreover,
%
\begin{equation}
\label{eqERW}
\rwE_N^\varepsilon[L(N)] =
\cases{
c_{47} \varepsilon N \bigl(1+o(1)\bigr) &\quad$(d=1)$,\vspace*{2pt}\cr
e^{-c_{48}/\varepsilon (1+o(1))} N &\quad$(d=2)$.}
\end{equation}
\end{theorem}
\begin{pf}
Notice first that in terms of the variables $\hat\tau_i$,
\[
\rwP_N^\varepsilon\Biggl( \sum_{i=1}^K \tau_i = N \Biggr) = \frac{\trwP
( \sum_{i=1}^K \hat\tau_i = N )}{\trwP(\exists K\geq1
\dvtx
\sum_{i=1}^K \hat\tau_i = N )} .
\]
By a standard large deviation estimate,
\[
\trwP\Biggl( \sum_{i=1}^K \hat\tau_i = N \Biggr) \leq e^{-c_{49}
(K\vee N)}
\]
for all $K$ such that $|K-N/\trwE[\hat\tau_1]| > \eta N$. Since
\[
\trwP\Biggl(\exists K\geq1 \dvtx \sum_{i=1}^K \hat\tau_i = N
\Biggr)\to
1/\trwE[\hat\tau_1] ,
\]
it thus follows that
\[
\rwP_N^\varepsilon\biggl( \biggl|\frac{L(N)}N - \frac1{\trwE[\hat
\tau
_1]}\biggr| \geq\eta\biggr) \leq e^{-c_{50} N}.
\]
\upqed
\end{pf}
%
\begin{cor}\label{corRWpinning}
Assume that $d=1$ or $d=2$. Then there exist $c_{51}, c_{52}>0$ such
that, for any small enough $\delta>0$, and $N$ large enough,
\[
{\rwP}\bigl(L(N) \geq\delta N\bigr) \geq
\cases{
e^{-c_{51} \delta^2 N}, &\quad if $d=1$,\cr
e^{-c_{52} (\delta/|{\log\delta}|) N}, &\quad if $d=2$.}
\]
\end{cor}
\begin{pf}
Using a well-known inequality (\cite{Giacomin}, (A.13))
\[
{\rwP}\bigl(L(N) \geq\delta N\bigr) \geq\rwP_N^\varepsilon\bigl(L(N) \geq\delta N\bigr)
\exp\biggl\{- \frac{H(\rwP_N^\varepsilon| {\rwP}) + e^{-1}} {\rwP
_N^\varepsilon(L(N) \geq\delta N)} \biggr\},
\]
where $H(\rwP_N^\varepsilon| {\rwP})$
denotes the relative entropy of $\rwP_N^\varepsilon$ w.r.t. ${\rwP}$.
We choose
\[
\varepsilon=\varepsilon(\delta) =
\cases{
c \delta &\quad$(d=1)$,\cr
c/|{\log\delta}| &\quad$(d=2)$,}
\]
with $c$ chosen in such a way that [remember (\ref{eqERW})]
\[
\rwE_N^\varepsilon[L(N)] \in(2\delta N,3\delta N).
\]
It then follows from (\ref{eqLLNRW}) that
\[
\rwP_N^\varepsilon\bigl(L(N) \geq\delta N\bigr) \geq\tfrac12
\]
for all $N$ large enough. But for large enough $N$,
\[
H(\rwP_N^\varepsilon| {\rwP}) = \varepsilon\rwE_N^\varepsilon[L(N)] -
\log
\rwZ_N^\varepsilon+\log\rwP(X_N=0)\leq3 \varepsilon\delta N.
\]
The conclusion follows.
\end{pf}
\end{appendix}

\section*{Acknowledgments}

S. Friedli gratefully acknowledges
the Section de\break Math\'{e}matiques of the University of Geneva for
hospitality while completing this project.



\printaddresses

\end{document}